\documentclass[10pt,onecolumn,letterpaper]{article}
\usepackage{cvpr}
\usepackage{times}
\usepackage{subfigure}
\usepackage{epsfig}
\usepackage{graphicx}
\usepackage{amsthm}
\usepackage{amsmath}
\usepackage{amssymb}
\usepackage{algorithm,algpseudocode}
\usepackage[bookmarks=ture]{hyperref}
\hypersetup{
	breaklinks	=true,
	bookmarks	=true,
	hyperindex=true, 
	colorlinks=true, 
	linkcolor=blue,
	citecolor=red, 
	urlcolor=red	
}
\usepackage{multirow}
\newtheorem{theorem}{Theorem}
\newtheorem{lemma}{Lemma}
\newtheorem{definition}{Definition}
\newtheorem{assumption}{Assumption}

\usepackage{datetime}
\usepackage{fancyhdr}
\usepackage{authblk}
\usepackage{fancyhdr}
\usepackage{pdftexcmds}
\def\r1{\color{black}}
\def\b1{\color{black}}
\def\B1{\bm}
\cvprfinalcopy % *** Uncomment this line for the final submission
\def\cvprPaperID{****} % *** Enter the CVPR Paper ID here
\allowdisplaybreaks

\begin{document}
	\lhead{}
	\lfoot{\date{\today},\date{\currenttime}}
	%\rfoot{SARC}
	
	\title{ Stochastic Second-order Methods for Non-convex Optimization with 
		Inexact Hessian and 
		Gradient}
	%econd-order method for convex and 
	%	non-convex function}
	\author[1]{ Liu Liu\thanks{lliu8101@uni.sydney.edu.au}}
	\author[2]{ Xuanqing Liu\thanks{xqliu@ucdavis.edu}}
	\author[2]{ Cho-Jui Hsieh\thanks{chohsieh@ucdavis.edu}}
	\author[1]{ Dacheng Tao\thanks{dacheng.tao@sydney.edu.au}}
	\affil[1]{UBTECH Sydney AI Centre and SIT, FEIT, The University of Sydney}
	\affil[2]{University of California, Davis}
	\maketitle
	
	\thispagestyle{fancy}
	\begin{abstract}
		Trust region and cubic regularization methods have demonstrated good performance	in 
		small 
		scale non-convex optimization, showing the ability to escape from saddle points. Each 
		iteration 
		of these methods involves computation of 	gradient, Hessian and function value in order 
		to 
		obtain search direction  and adjust the radius or cubic regularization parameter. However, 
		exactly computing those quantities are 	too expensive in large-scale problems such as 
		training 
		deep networks. In this paper, we study a family of stochastic trust region and cubic 	
		regularization methods when gradient, Hessian and 	function values are computed inexactly, 
		and 
		show the iteration complexity	to achieve $\epsilon$-approximate second-order optimality 
		is in 
		the same order with previous work for which gradient and function values are computed 
		exactly. 
		The mild conditions on inexactness can be achieved in finite-sum minimization using random 
		sampling. We show the algorithm performs well on training convolutional neural networks 
		compared with previous second-order methods.
	\end{abstract}
	%%%%%%%%%%%%%%%%%%%%%%%%%%%%%%%%%%%%%%%%%%%%%%%%%%%%%%%%%%%%%%%%%%%%%%%%%%%%%%%%%%%%%%%%%%%%%
	%%%%%%%%%%%%%%%%%%%%%%%%%%%%%%%%%%%%%%%%%%%%%%%%%%%%%%%%%%%%%%%%%%%%%%%%%%%%%%%%%%%%%%%%%%%%%%%
	%%%%%%%%%%%%%%%%%%%%%%%%%%%%%%%%%%%%%%%%%%%%%%%%%%%%%%%%%%%%%%%%%%%%%%%%%%%%%%%%%%%%%%%%%%%%%%%%
	\section{introduction}
	In this paper, we consider the unconstrained optimization problem:
	\begin{align}\label{Newton:Problem}
	\mathop {\min }\limits_{x \in {\mathbb R^d}} f( x ) = 
	\frac{1}{n}\sum\limits_{i = 1}^n {{f_i}( x )},
	\end{align}
	where $f_i(x)$ is smooth and not necessarily convex. Such finite-sum
	structure is increasingly popular in modern machine learning tasks, 
	especially in deep learning, where each $f_i(\cdot)$ corresponds to the loss 
	of a training sample. For large-scale problems, computing the full gradient and Hessian is 
	prohibitive,  so Stochastic gradient descent (SGD) has 
	become the most popular method. 
	%become an important 
	%first-order optimization method. 
	%and its convergence rate has been well studied in both
	%convex and non-convex cases \CH{cite some SGD papers}.
	First order methods such as gradient descent and SGD are guaranteed to 
	converge to stationary points, which can be a saddle point or a local 
	minimum\footnote{Some recent analysis indicates SGD can escape from saddle 
		points in certain cases~\cite{ge2017no,jin2017escape}, but SGD is not the focus of this 
		paper. }. 
	%  Very recently, second order methods have been applied to deep learning, and some of them 
	%  has demonstrated better performance in terms of number of epochs, but the convergence 
	%properties
	%  are y not well established. 
	%  
	%  In order to scale to large datasets, 
	%  both Hessian and gradient needs to be subsampled. However, the theoretical guarantee
	%  for second order methods under inexact gradient and Hessian is not well studied. Existing 
	%analysis
	%  either has strong assumption on Hessian \CH{cite} or assume there is no noise in gradient 
	%computation \CH{cite}. 
	
	It is known that second-order methods, by utilizing the Hessian information, 
	can more easily escape from saddle points.
	At each iteration, second-order methods typically build a quadratic 
	approximation function around the current solution $x_k$ by 
	\begin{align}\label{Newton:Definition_m}
	{m_k}( s) = f( {{x_k}} ) + \langle 
	{{g(x_k)},s} \rangle  + \frac{1}{2}\langle {s,{B(x_k)}s} 
	\rangle,
	\end{align}
	where $g(x_k)$ is the approximated gradient and $B(x_k)$ is 
	the symmetric matrix. To update the current solution, a common strategy is to 
	minimize this quadratic approximation within a small region. Algorithms based 
	on this idea including Trust Region method (TR)~\cite{conn2000trust} and 
	Adaptive Regularization using Cubics 
	(ARC)~\cite{cartis2011adaptivea,cartis2011adaptiveb} 
	have demonstrated good performance on small-scale non-convex problems.
	
	However, for large-scale optimization such as training deep neural networks, it is impossible 
	to 
	compute gradient and Hessian {\it exactly} for every update. As a result, stochastic 
	second-order 
	methods have been studied in the past few years.  \cite{xu2017newton, xu2017second} proposed
	TR and ARC methods with inexact Hessian, in which the second-order information is approximated 
	by  the subsampled Hessian matrix, yet the gradient is still computed exactly. 
	%Moreover, they also require the objective function to be computed 
	%at each step in order to adjust the trust region radius.
	%In terms of cubic regularization method, 
	%The second problem is that updating the radius depend on computing the 
	%objective value in (\ref{SFSNewton:Problem1}), which greatly affect the 
	%efficiency of computation,  especially in deep learning.
	\cite{kohler2017sub} proposed a stochastic version of ARC, but they require a much stronger 
	condition in both gradient and Hessian approximation, thus they need to keep increasing the 
	sample size as iteration goes. 
	%needs to increase as iteration goes and it will eventually reduce
	%to the full gradient/Hessian version. 
	More recently, \cite{tripuraneni2017stochastic} provide 	stochastic cubic regularization 
	method, but they do not have an adaptive way 	to adjust the regularization parameters.
	%even though avoiding the computation 
	%of function, the parameters in cubic regularization term is not known in 
	%advance.
	
	In this paper, we consider a simple and practical stochastic version of trust region 
	and cubic regularization methods (denoted as STR and SARC respectively). In the quadratic 
	approximation~\eqref{Newton:Definition_m}, 
	we replace both $g(x_k)$ and $B(x_k)$ by the approximate gradient and Hessian
	with a fixed approximation error, which can be achieved using a fixed sample 
	size. 
	Furthermore, we consider the most general case where 
	the trust region radius or cubic regularization parameter is adaptively adjusted
	by checking the {\it subsampled} objective function value. 
	
	Note that we do 
	not claim STR and SARC are ``new'' algorithms, since it is natural to transform TR 
	and ARC to the stochastic setting. The question is whether this simple 
	idea works in theory and in practice, and our contribution is to provide an affirmative 
	answer to this question. 
	Our contribution 
	can be summarized as follows: 

	\begin{itemize}
		\item We provide a theoretical analysis of convergence and iteration 
		complexity for STR and SARC. Unlike~\cite{kohler2017sub}, we do not require
		the approximation error of Hessian and gradient estimation to be related 
		to $s$ (the update step). Furthermore, even though the proof framework is similar to
		\cite{xu2017second}, we provide novel analysis to model the case when 
		both gradient and Hessian are inexact, while \cite{xu2017second} does 
		not allow inexact gradient. 
		%		framework is based on \cite{xu2017second}, that separately consider 
		%		three situations related to gradient and Hessian, our proof method for 
		%		escaping the saddle point in the non-convex problem is different 
		%		\cite{xu2017second}. Because we have more conditions for sub-sampling 
		%		the gradient, which will reduce the cost of gradient computation.
		
		\item We use the operator-Bernstein inequality to bound the sub-sampled 
		function value, 
		which 
		enables 
		automatically adjusting trust region radius and cubic regularization parameter 
		using subsamples. Even when function value, gradient  and Hessian
		are all inexact, we are able to show that the iteration complexity is in the same order
		with \cite{xu2017second,kohler2017sub}. 
		%      that apply 
		%		operator-Bernstein inequality to obtain the function value. Though the 
		%		bound of conditions is related to the radius or the parameter of cubic 
		%		regularization $\sigma$, we have proof that there exit 	minimal radius 
		%		and maximal parameter $\sigma$ such that sub-sampling the function with 
		%		fix size under probability. Furthermore, we also give the theoretical 
		%		analysis, which  preserves the iteration complexity as in 
		%		\cite{xu2017second, 
		%		kohler2017sub}.
		\item We conduct experiments on CIFAR-10 data with VGG network and show that
		the proposed algorithms are faster than  the existing trust region and cubic 
		regularization methods in terms of running time.  
	\end{itemize}
	\subsection{Our results}
	We present the iteration complexity for both proposed methods under the  Assumptions in 
	Preliminary where the parameters are defined as well.
	
	\textbf{{STR}} The total number of iterations is 
	$\mathcal{O}( \max \{ (\varepsilon _{\nabla f} - \varepsilon _g)^{ - 
		2},(\varepsilon _H - \varepsilon _B)^{ - 3} \} )$\footnote{We use 
		$\mathcal{O}(\cdot)$ to hide constant factors.}.
	
	\textbf{{SARC}} The iteration complexity is the same order of  STR.
	If the condition of the terminal criterion is satisfied, then the total number of iterations is 
	\[\mathcal{O}( \max \{ (\varepsilon _{\nabla f} - \varepsilon _g)^{ - 
		3/2},(\varepsilon _H - \varepsilon _B)^{ - 3}\}  ),\]
	which is better than STR.
	%%%%%%%%%%%%%%%%%%%%%%%%%%%%%%%%%%%%%%%%%%%%%%%%%%%%%%%%%%%%%%%%%%%%%%%%%%%%%%%%%%%%%%%
	%%%%%%%%%%%%%%%%%%%%%%%%%%%%%%%%%%%%%%%%%%%%%%%%%%%%%%%%%%%%%%%%%%%%%%%%%%%%%%%%%%%%%%%
	%%%%%%%%%%%%%%%%%%%%%%%%%%%%%%%%%%%%%%%%%%%%%%%%%%%%%%%%%%%%%%%%%%%%%%%%%%%%%%%%%%%%%%%
	\subsection{Related Work}
	With the increasing size of data and model, stochastic optimization becomes more and more 
	popular since computing the gradient and Hessian are prohibitively 
	expensive. 
	%Especially the burst development of deep learning and its 
	%application, many researchers turn their focus on theoretical analysis and 
	%try to apply the classical optimization methods into stochastic versions. 
	
	For the stochastic first-order optimization, stochastic gradient descent (SGD)
	\cite{zhang2004solving, shalev2011pegasos, ghadimi2016accelerated} is 
	absolutely the main 
	method especially in training deep neural networks and other large-scale machine learning 
	problems, due to its simplicity and effectiveness. 
	%implement and unbiased estimation of the 
	%gradient. 
	However, the estimated gradient will induce the noise such that the 
	variance of the gradient may not approximate zero even when converged
	to a stationary point. 
	Stochastic variance reduction gradient (SVRG) 
	\cite{johnson2013accelerating} and SGAG \cite{defazio2014saga} are two 
	typical methods to reduce the variance of the gradient estimator, which lead to faster
	convergence especially in the convex setting.  
	%that divide the iteration into epochs, in which 
	%the full gradient and one or block stochastic gradient are added to the 
	%gradient $g(x)$ in (\ref{SFSNewton:Definition_h}) so that the gradient estimator
	%remains unbiased and the variance also gradually reduces to 0. 
	Several other related variance reduction methods are developed and analyzed 
	for non-covex problems
	%Ever since the related variance reduction papers are emerging such as 
	\cite{reddi2016stochastic,allen2016variance}. However, 
	using only first-order information, the saddle point may not be escaped 
	even though \cite{anandkumar2016efficient} prove that SGD with noise can 
	escape but under certain conditions. 
	%strict saddle function.

	For second-order optimization, Newton-typed methods rely on building a 
	quadratic approximation around the current solution, and 
	by exploring the curvature information it can better avoid saddle points in 
	non-convex optimization. 
	%stochastic Newton method is often 
	%used to the quadratic model. Explore the curve information is important in 
	%the  
	%non-convex problem, especially near saddle point. 
	The negative eigenvector 
	of the Hessian information provides the decrease direction for the updates.
	Using exact Hessian is often time consuming, so 
	Broyden-Fletcher-Goldfarb-Shanno (BFGS) and 
	Limited-BFGS \cite{nocedal1980updating} are two widely methods that approximate
	Hessian using first-order information. 
	Another important technique for Hessian approximation is 
	sub-sampling the function $f_i(\cdot)$ to obtain the estimated 
	Hessian. Both \cite{byrd2011use} and \cite{erdogdu2015convergence} using 
	the stochastic Hessian matrix to obtain the global convergence, while the 
	former requires $f_i$ to be smooth and strongly convex.
	Furthermore, \cite{roosta2016suba,roosta2016subb} apply the 
	sub-sampled gradient and Hessian to the quadratic model and give the 
	convergence analysis of second-order methods thoroughly and 
	quantitatively. 
	%However, how to determine whether the approximated 
	%quadratic function is good or not, will lead us to consider another 
	%class method that through adding constraint condition or regularization 
	%term.
	
	Trust region Newton method is a classical second-order method that 
	searches the update direction only within a trust region around the current 
	point. 
	The size of trust region is critical to the effectiveness of 
	search direction. The region will be updated based on measuring whether
	the quadratic approximation is an adequate representation of the function or not. 
	Following the sub-sampling for Hessian matrix as in 
	\cite{roosta2016suba,roosta2016subb}, 
	\cite{xu2017newton, xu2017second} apply such inexact Hessian to trust 
	region 	method and also provide the convergence and iteration complexity. 
	Similar to the trust region method,	
	\cite{cartis2011adaptivea,cartis2011adaptiveb} 
	introduced adaptive 
	cubic regularization methods for unconstrained optimization, in which the 
	Hessian metrics can be replaced by the approximate matrix. 
	\cite{kohler2017sub} apply the operator-Bernstein inequality to approximate 
	the Hessian matrix and gradient into a quadratic function with cubic 
	regularization.	However, the sample approximate condition is subject to the 
	search direction $s$, thus they need to increase the sample size at each step. 
	To overcome this issue, \cite{xu2017newton} provide 
	another approximation condition of Hessian matrix that does not depend on the 
	search step $s$. However, they assume gradient has to be computed exactly, 
	which is not feasible in large-scale applications. 
	Furthermore, each update of the $\rho$ that measure the adequacy  
	of the function  will need the full computation of the objective function, 
	which will lead to more computation cost. 
	%{Recently, Nilesh \cite{tripuraneni2017stochastic} introduce the 
	%stochastic 
	%	cubic regularization for the non-convex problem.}

	The rest of paper is organized as follows. Section  
	\ref{Newton:Section:Preliminary} 
	gives the  preliminary about the assumptions and definition. The 
	sub-sampling method for estimating the corresponding function, gradient and 
	Hessian is in Section \ref{Newton:Section:subsampling}.
	Section \ref{Newton:Section:STR} and  \ref{Newton:Section:SARC} 
	respectively present the stochastic trust region method and cubic 
	regularization method, and their convergence and iteration complexity.  
	Section 
	\ref{Newton:Section:Experiment} gives the experimental results. Section 
	\ref{Newton:Section:Conclusion} concludes our paper.
	%%%%%%%%%%%%%%%%%%%%%%%%%%%%%%%%%%%%%%%%%%%%%%%%%%%%%%%%%%%%%%%%%%%%%%%%%
	%%%%%%%%%%%%%%%%%%%%%%%%%%%%%%%%%%%%%%%%%%%%%%%%%%%%%%%%%%%%%%%%%%%%%%%%
	%%%%%%%%%%%%%%%%%%%%%%%%%%%%%%%%%%%%%%%%%%%%%%%%%%%%%%%%%%%%%%%%%%%%%%%%
	%%%%%%%%%%%%%%%%%%%%%%%%%%%%%%%%%%%%%%%%%%%%%%%%%%%%%%%%%%%%%%%%%%%%%%%%%
	\section{Preliminary}\label{Newton:Section:Preliminary}
	%Through this paper, we denote by ${\nabla f( x )}$ the gradient of function 
	%$f(x)$, and ${{\nabla ^2}f( x )}$ the Hessian of function $f(x)$. 
	For a 
	vector $x$ and a matrix $X$,  we use $\|x\|$ and $\|X\|$ to denote the 
	Euclidean norm and the matrix spectral norm, respectively. We use $\mathcal{S}$ to denote
	the set and $|\mathcal{S}|$ to denote its cardinality. 
	For the matrix $X$, we use $\lambda_{\text{min}}(X)$ and $\lambda_{\text{max} 
	}(X)$ to denote its smallest and largest eigenvalue.  
	In the following, we give  assumptions and definition about the characteristic of function, the 
	approximate conditions, related bounds, and optimality definition.
	
	\begin{assumption}(Lipschitz Continuous)\label{Newton:Assumption:function} For the function 
		$f(x)$, we assume that ${\nabla ^2}f( x )$ and ${\nabla}f( x )$  are 
		Lipschitz continuous 
		satisfying $\| \nabla ^2f( x ) - {\nabla ^2}f( y ) \| \le L_H\| x - y\|$ 
		and $\| \nabla 
		f( x ) - {\nabla}f( y ) \| \le L_{\nabla f}\| x - y\|$, $\forall x,y \in 
		{\mathbb{R}^d}$.
		%	\begin{align*}
		%%%%%	\label{SFSNewton:Assumption:function:FunctionLipschitz}
		%%%%	\| { f( x ) -  f( y )} 	\| 	&\le {L_f}\| {x - y} \|,\\
		%%%%%	\label{SFSNewton:Assumption:function:GradientLipschitz}
		%%%%	\| {\nabla f( x ) - \nabla f( y )} 	\| 	&\le {L_{\nabla f}}\| {x - y} \|,\\
		%%%%%	\label{SFSNewton:Assumption:function:HessianLipschitz}
		%%%%	\| {{\nabla ^2}f( x ) - {\nabla ^2}f( y )} \| &\le L_H\| {x - y} \|.
		%	\end{align*}
	\end{assumption}
	
	\begin{assumption}(Approximate)\label{Newton:Assumption:approximation} For function $f(x)$, 
	the  
		approximate gradient $g(x)$ and Hessian matrix $B(x)$ satisfy
		\begin{align}\label{Newton:Assumption:approximation-gB}
		\left\| {\nabla f(x) - g(x)} \right\| \le {\varepsilon _g},
		\left\| {{\nabla ^2}f(x) - B(x)} \right\| \le {\varepsilon _B}.
		\end{align}	
		with ${\varepsilon _g},{\varepsilon _B} > 0$. The approximated function $h(x)$ at 
		$k$-iteration 
		satisfies 
		\begin{align}\label{Newton:Assumption:approximation-h}
		\r1 \left\| {f(x_k) - h(x_k)} \right\| \le {\varepsilon _h}{\|s_k\|^2},{\varepsilon _h} > 0.
		%\left\| {f(x_k) - h(x_)} \right\| \le {\varepsilon _h}/{\sigma _k}.
		\end{align}
		%with ${\varepsilon _h} > 0$, the parameters ${\Delta _k}$ and ${\sigma _k}$	are defined 
		%in 
		%Algorithm \ref{Newton:STR:Algorithm} and \ref{Newton:SARC:Algorithm}.
	\end{assumption}
	
	\begin{assumption}(Bound)\label{Newton:Assumption:bound} 
		For $i\in[n]$, the bound assumptions are the 
		function $f_i(x)$ satisfies $\left\| 
		{{f_i}({x})} \right\| \le {\kappa _f}$, $\left\| {\nabla {f_i}({x})} \right\| \le {\kappa 
			_{\nabla f}}$, and $\left\| {{\nabla ^2}{f_i}({x})} \right\| \le {\kappa _H}$.
	\end{assumption}
	
\begin{figure*}[t]
		\centering
		\includegraphics[width=0.42\textwidth]{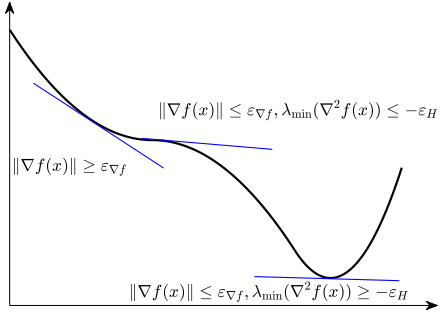}
		\caption{Illustration of three situations in analyzing the convergence 
			to first and second critical 
			point.}\label{Newton:Figure3Situation}
	\end{figure*}

	\begin{assumption} (Bound)\label{Newton:Assumption:Bound:Variance}
		We assume that $H_1$ and $H_2$ are the upper bounds on the variance of the $\nabla f_i(x)$ 
		and 
		$\nabla^2f_i(x)$, that is	
		\begin{center}
			$
			\frac{1}{n}\sum\limits_{i = 1}^n {{{\left\| {\nabla {f_i}\left( x \right) - \nabla 
			f\left( 
							x \right)} \right\|}^2}}  \le {H_1^2},\frac{1}{n}\sum\limits_{i = 1}^n 
			{{{\left\| 
						{{\nabla 
								^2}{f_i}\left( x \right) - {\nabla ^2}f\left( x \right)} 
								\right\|}^2}}  \le 
			{H_2^2}.
			$
		\end{center}
	\end{assumption}
	The following lemmas are important tools for analyzing the convergence of the proposed 
	algorithm, 
	which are used to characterize the variance of random variable decreasing with the factor 
	related 
	to the set size.
	%%%==============================================================================
	%%%===============================================================================
	%%%=============================================================================
	\begin{lemma}\label{Newton:Appendix:Tool:RandomSubset}
		If $v_1,...,v_n\in \mathbb{R}^d$ satisfy 
		$\sum\nolimits_{i = 1}^n {{v_i}}  = \vec 0$, and $\cal A$ is a non-empty, 
		uniform 
		random subset of $[n]$, $A=|\mathcal{A}|$, then
		\begin{center}
			${\mathbb{E}_{\mathcal{A}}} {{{\left\| {\frac{1}{A}\sum\nolimits_{b \in 
								{\mathcal{A}}} 
							{{v_b}} 
						} 
						\right\|}^2}}  \le \frac{{\mathbb{I}\left( {A < n} 
					\right)}}{A}\frac{1}{n}\sum\limits_{i = 1}^n {v_i^2}.$
		\end{center}
	\end{lemma}
	\begin{definition}(($\varepsilon_{{\nabla f}},\varepsilon _H$)-Optimality). 
		Given 
		${\varepsilon _{{\nabla f}}},{\varepsilon _H} \in [0,1]$, $x$ is an $( 
		{{\varepsilon_{{\nabla f}}},{\varepsilon _H}}))$-Optimality solution to 
		problem 
		(\ref{Newton:Problem}), if
		\begin{align*}
		\| 	{\nabla f( {{x}} )} \| \le 	{\varepsilon _{\nabla 
				f}},\,\lambda_\text{min} ({\nabla ^2}f( x 
		) 
		) 	\ge  - 	{\varepsilon_H }.
		\end{align*}
	\end{definition}
	Furthermore, we introduce three index sets
	\begin{align*}
	%\label{SFSNewton:SARC:Definition:S_g}
	\mathcal{S}_{{\nabla f}}&\mathop  = \limits^\text{def}\{x: \,  \|{\nabla 
		f( {{x}} )} \| \ge {\varepsilon _{\nabla f}}\}, \\
	%\label{SFSNewton:SARC:Definition:S_H}	
	\mathcal{S}_{H}&\mathop  = \limits^\text{def}\{x: \|{\nabla 
		f( {{x}} )} \| \ge {\varepsilon _{\nabla f}} \,\text{and} 
	\,\lambda_\text{min} 
	({\nabla^2f(x )} ) \le - {\varepsilon_H }\}, \\
	%\label{SFSNewton:SARC:Definition:S_star}
	\mathcal{S}_{*}&\mathop  = \limits^\text{def} \{x:\, \|
	\nabla f( {{x}} ) \| \le \varepsilon_{\nabla f}\,
	\text{and}\, \lambda_\text{min} (\nabla^2f( x )) \ge  -\varepsilon_H\}.	
	\end{align*}
	where  {$ {\varepsilon _{\nabla f}}>\varepsilon_g$ and 
		$\varepsilon_H>\varepsilon_B$}. In order to clearly classify three 
	situations, we give a simple 
	geometry 
	illustration, as shown 
	in Figure \ref{Newton:Figure3Situation}.
	%%%%%%%%%%%%%%%%%%%%%%%%%%%%%%%%%%%%%%%%%%%%%%%%%%%%%%%%%%%%%%%%%%%%%%%%%
	%%%%%%%%%%%%%%%%%%%%%%%%%%%%%%%%%%%%%%%%%%%%%%%%%%%%%%%%%%%%%%%%%%%%%%%%
	%%%%%%%%%%%%%%%%%%%%%%%%%%%%%%%%%%%%%%%%%%%%%%%%%%%%%%%%%%%%%%%%%%%%%%%%
	%%%%%%%%%%%%%%%%%%%%%%%%%%%%%%%%%%%%%%%%%%%%%%%%%%%%%%%%%%%%%%%%%%%%%%%%%
	\section{Sub-sampling for finite-sum 
		minimization}\label{Newton:Section:subsampling}
	For the finite-sum problem~\eqref{Newton:Problem}, 
	we can estimate 
	$f(x)$, $\nabla f(x)$, and $\nabla^2 f(x)$ by random sub-sampling,  
	%a sample of index from $[n]$ to form the sampling set. 
	which can drastically reduce the computational complexity.
	Here, we use $\mathcal{S}_h$, $\mathcal{S}_g$ and $\mathcal{S}_B$ to denote the sample 
	collections 
	for estimating $f(x)$, $\nabla f(x)$ and $\nabla^2f(x)$, respectively, where $\mathcal{S}_h$, 
	$\mathcal{S}_g$ and 
	$\mathcal{S}_B\subseteq [n]$. The approximated functions are formed by
	\begin{align}
	\label{Newton:STR:Subsampling:h}
	h(x)&=\frac{1}{\left| {{\mathcal{S}_h}} \right|}\sum\nolimits_{i \in {\mathcal{S}_h}} 
	{{f_i}\left( 
		x \right)},\\
	\label{Newton:STR:Subsampling:g}
	g( {{x}} ) &= \frac{1}{{\left| {{\mathcal{S}_g}} \right|}}\sum\nolimits_{i \in 
		{\mathcal{S}_g}}^{} {\nabla 
		{f_i}\left( 
		{{x}} \right)},\\
	\label{Newton:STR:Subsampling:B}
	B( {{x}} ) &= \frac{1}{{| {{\mathcal{S}_B}} |}}\sum\nolimits_{i \in {\mathcal{S}_B}}^{} 
	{{\nabla 
			^2}{f_i}( {{x}} )}.
	\end{align}
	Most papers use operator-Bernstein inequality to probabilistically guarantee such properties, 
	such 
	as \cite{drineas2006fast,mahoney2011randomized} use  approximate matrix multiplication results 
	as a 
	fundamental primitive 
	in RandNLA  \cite{drineas2006fast,mahoney2011randomized} to control the 
	approximation error of $\nabla f(x)$.  Furthermore,  the vector-Bernstein inequality 
	\cite{candes2011probabilistic,gross2011recovering} is applied in \cite{ kohler2017sub} to 
	obtain 
	the of sub-sample bound of the gradient, which is different from 	\cite{drineas2006fast}. 
	However, 
	the number of sub-samples depend on the search direction $s$ in advance, which will affect the 
	estimation of sub-sampling. 	Replacing the condition by  
	(\ref{Newton:Assumption:approximation-gB}), we can have
	%%%==============================================================================
	%%%===============================================================================
	%%%=============================================================================
	\begin{lemma} 
		Suppose the Assumption \ref{Newton:Assumption:bound} holds,
		if $| {{\mathcal{S}_g}} | \ge 16\log ( 
		{{{2d}}/{\delta }} 
		){{L_f^2}}/{{\varepsilon	_g^2}}$, then $g(x)$ formed by 
		(\ref{Newton:STR:Subsampling:g})  
		satisfies $\| \nabla f( x ) - g(x) \| \le {\varepsilon}_{g}$ with probability 
		($1-\delta$).
	\end{lemma}
	Note that, we do not give the proof of above Lemmas as the difference lies on different 
	conditions. 
	We also obtain the approximate gradient $g(x)$ based on above results. 
	%Using operator-Bernstein inequality in \cite{gross2010note},  \cite{xu2017newton} provide the 
	%sub-sample bound about the approximated Hessian. 
	%Since the condition for approximation in (\ref{Newton:Assumption:approximation-gB}), 
	%we use the sum-sample bound directly that can be applied into our algorithms.
	%%%==============================================================================
	%%%===============================================================================
	%%%=============================================================================
	\begin{lemma}
		Suppose the Assumption \ref{Newton:Assumption:bound} holds,
		If $| \mathcal{S}_B | \ge \log ( {{{2d}}/{\delta }} ) 
		{{16L_B^2}}/{{\varepsilon_B^2}}$, then  $B(x)$ formed by (\ref{Newton:STR:Subsampling:B})
		satisfies $\| {\nabla ^2}f( x ) -{B(x)}\| \le 	{\varepsilon}_B$ with probability 
		($1-\delta$).
	\end{lemma}
	In order to reduce the computation cost for updating the $\rho$, we also use the sub-sampling 
	combing with operator-Bernstein inequality to obtain the approximate function probabilistically 
	satisfying (\ref{Newton:Assumption:approximation-h}). Note that, before obtain the approximate 
	function $f(x)$, $s$ has been solved. Thus, we can use $\|s\|$ directly. Furthermore, if the 
	approximate function $h(x)$ satisfies the condition in 
	(\ref{Newton:Assumption:approximation-h}), 
	but could not be guaranteed the bound of $\left\| {\nabla h\left( x \right) - \nabla f\left( x 
		\right)} 
	\right\|$ and $\| {{\nabla ^2}h( x ) - {\nabla ^2}f( x )} \|$, 
	which will be later used to analyze the radius. Thus, we present a important assumption 
	and use Lemma \ref{Newton:Appendix:Tool:RandomSubset} to derive the upper bound. 
	%%%==============================================================================
	%%%===============================================================================
	%%%=============================================================================
	\begin{lemma}\label{Newton:Subsample:Lemma} 
		Suppose the Assumption \ref{Newton:Assumption:bound} and 
		\ref{Newton:Assumption:Bound:Variance} hold,
		If 	 $|{{\cal S}_h}| \ge \log (2d/\delta )16\kappa _f^2/( {\varepsilon _h^2{{\| {{s_k}} 
		\|}^4}} 
		)$,
		then 
		$h(x)$	formed by (\ref{Newton:STR:Subsampling:h}) satisfies $\|  f(x_k ) - h(x_k) \| \le 
		{\varepsilon}{_h}\|s_k\|^2$ with probability 
		($1-\delta$).
		Furthermore,  we can also have the upper bounds with the gradient and the Hessian of $h(x)$,
		\[ 	{\left\| {\nabla h\left( x \right) - \nabla f\left( x \right)} \right\|^2} \le 
		\frac{{\mathbb{I}\left( {\left| {{\mathcal{S}_h}} \right| < n} \right)}}{{\left| 
				{{\mathcal{S}_h}} \right|}}{H_1},
		{\left\| {{\nabla ^2}h\left( x \right) - {\nabla ^2}f\left( x \right)} \right\|^2} \le
		\frac{{\mathbb{I}\left( {\left| {{\mathcal{S}_h}} \right| < n} \right)}}{{\left| 
				{{\mathcal{S}_h}} \right|}}{H_2}.\]
	\end{lemma}
	%%%%%%%%%%%%%%%%%%%%%%%%%%%%%%%%%%%%%%%%%%%%%%%%%%%%%%%%%%%%%%%%%%%%%%%%%%
	%%%%%%%%%%%%%%%%%%%%%%%%%%%%%%%%%%%%%%%%%%%%%%%%%%%%%%%%%%%%%%%%%%%%%%%%%%%
	%%%%%%%%%%%%%%%%%%%%%%%%%%%%%%%%%%%%%%%%%%%%%%%%%%%%%%%%%%%%%%%%%%%%%%%%%%%%
	\section{Stochastic Trust Region method}\label{Newton:Section:STR}
	\begin{algorithm*}[t]
		\caption{STR  with Inexact Hessian and Gradient}
		\label{Newton:STR:Algorithm}
		\begin{algorithmic}[1]
			\Require given $x_0$, $r_2\ge 1> r_1$, $1>\eta>0$, $\varepsilon>0$ and 
			$\Delta_0,\Delta_\text{max}>0$.
			\Ensure
			\For{$k$=1 to T }
			\State Compute the approximate gradient $g(x_k)$ and 
			Hessian matrix 	$B(x_k)$ based on 
			(\ref{Newton:STR:Subsampling:g}) 
			and (\ref{Newton:STR:Subsampling:B}).
			\State 
			\Comment({\r1 if 
				%${\lambda _{\text{min} }}( {\nabla^2f(x_k)} ) 	\le  - 	{\varepsilon _H}$, 
				%$||\nabla f\left( {{x_k}} \right)|| \le {\varepsilon _{\nabla f}}$, 
				$\|g\left( {{x_k}} \right)\| \le {\varepsilon 	_{\nabla 	f}}+ {\varepsilon _g}$, 
				we set $\mathcal{S}_h=\mathcal{S}_g$})
			\State Compute the direction vector $s_k$
			\begin{align}\label{Newton:STR:Algorithm-solution}
			s_k=\text{Subproblem-Solver}(g_k,B_k,\Delta_k).
			\end{align}	
			\State Compute the approximate  function $h(x_k)$ and $h(x_k+s_k)$ based on 
			(\ref{Newton:STR:Subsampling:h}).	
			\State 	Compute $
			{{\tilde \rho }_k} = \frac{{h\left( {{x_k}} \right) - h\left( {{x_k} + {s_k}} 
					\right)}}{{{m_k}\left( 0 \right) - {m_k}\left( {{s_k}} \right)}}
			$
			\State Set $\rho={{{\tilde \rho }_k} - \frac{{2{\varepsilon _f}{\r1\|s_k\|^2}}} 
				{{{m_k}\left( 0 \right) - {m_k}\left( {{s_k}} \right)}}}$			
			\State set ${x_{k + 1}} = \left\{ {\begin{array}{*{20}{l}}
				{{x_k} + {s_k},}&{\rho  \ge \eta ,}\\
				{{x_k},}&{{\rm{otherwise}}{\rm{.}}}
				\end{array}} \right.$
			\State set
			${\Delta _{k + 1}} = \left\{ {\begin{array}{*{20}{l}}
				{\text{min} \left\{ {{\Delta _{\text{max} }},{r_2}{\Delta _k}} \right\},}&{\rho  > 
				\eta 
					,}\\
				{{r_1}{\Delta _k},}&{{\rm{otherwise}}.}
				\end{array}} \right.$	
			\EndFor
			%\State \r1 how to terminal?
		\end{algorithmic}
	\end{algorithm*}
	
	In this section, we consider the stochastic trust region method for 
	solving the constrained optimization problem. At $k$-iteration, the 
	objective function is approximated by a quadratic model within a trust 
	region,
	\begin{align}\label{Newton:STR:Objective}
	\mathop {\min }\nolimits_{s\in \mathbb{R}^d} \,{m_k}( s 
	),\,\,\,\,	\text{subject to}\,\| s \| \le {\Delta _k},
	\end{align}
	where $\Delta_k$ is the radius, and ${m_k}( s )$ is defined in (\ref{Newton:Definition_m}). The 
	approximated gradient and Hessian 
	are formed based on the conditions in Assumption \ref{Newton:Assumption:approximation}. The 
	conditions of $\varepsilon_B $ in (\ref{Newton:Assumption:approximation-gB}) does not depend on 
	the 
	search direction $s$, which is the same as in 
	\cite{xu2017newton,xu2017second}. Moreover, we also define a new parameter $\varepsilon_g$, 
	which 
	has the same characteristic as $\varepsilon_B$.    In addition, the computations of $f(x)$ and 
	$f(x+s)$ are expensive as the objective function is finite-sum structure. Different from 
	\cite{xu2017newton} and \cite{kohler2017sub}, we  replace $f(x)$ and $f(x+s)$ with the 
	approximate 
	function $h(x)$ and $h(x+s)$ under the condition (\ref{Newton:Assumption:approximation-h}).  
	This 
	condition is  subjected to the search direction $s$. However, the  approximate 
	function $h(x)$ and $h(x+s)$ can be derived after obtaining the solution $s$ through 
	Subproblem-Solver.  Algorithm 
	\ref{Newton:STR:Algorithm} presents the process for updating the $x$ and $\Delta$.	This 
	section consists of two parts: Firstly, we analyze the role of radius $\Delta$ to ensure that 
	the 
	radius has the lower bound. Secondly, we derive the corresponding iteration complexity under 
	the 
	assumptions we present in Preliminary \ref{Newton:Section:Preliminary}.
	
	%%%%%%%%%%%%%%%%%%%%%%%%%%%%%%%%%%%%%%%%%%%%%%%%%%%%%%%%%%%%%%%%%%%%%%%%%%%%%%%%%%%%%%%
	%%%%%%%%%%%%%%%%%%%%%%%%%%%%%%%%%%%%%%%%%%%%%%%%%%%%%%%%%%%%%%%%%%%%%%%%%%%%%%%%%%%%%%%
	%%%%%%%%%%%%%%%%%%%%%%%%%%%%%%%%%%%%%%%%%%%%%%%%%%%%%%%%%%%%%%%%%%%%%%%%%%%%%%%%%%%%%%%
	\subsection{Bounds analysis  of radius }
	First of all, we present three important definition: $\rho_k$, $\tilde{\rho}$ and $\rho$. The 
	first two terms are defined as
	% \begin{center}
	% 	${\rho _k} = \frac{{f\left( {{x_k}} \right) - f\left( {{x_k} + {s_k}} 
	%\right)}}{{{m_k}\left( 0 
	% 			\right) - {m_k}\left( {{s_k}} \right)}},{{\tilde \rho }_k} = \frac{{h\left( {{x_k}} 
	% 			\right) 
	% 			- 
	% 			h\left( {{x_k} + {s_k}} \right)}}{{{m_k}\left( 0 \right) - {m_k}\left( {{s_k}} 
	% 			\right)}}.$
	% \end{center}
	\begin{align*}
	{\rho _k} = \frac{{f\left( {{x_k}} \right) - f\left( {{x_k} + {s_k}} \right)}}{{{m_k}\left( 0 
			\right) - {m_k}\left( {{s_k}} \right)}},{{\tilde \rho }_k} = \frac{{h\left( {{x_k}} 
			\right) 
			- 
			h\left( {{x_k} + {s_k}} \right)}}{{{m_k}\left( 0 \right) - {m_k}\left( {{s_k}} 
			\right)}}.
	\end{align*}
	Based on the inequality (\ref{Newton:Assumption:approximation-h}) in Assumption 
	\ref{Newton:Assumption:approximation}, we have
	%\begin{center}
	%	$\left| {{\rho _k} - {{\tilde \rho }_k}} \right| =\left| {\frac{{f\left( {{x_k}} \right) - 
	%				h\left( {{x_k}} \right) - \left( {f\left( {{x_k} + {s_k}} \right) - h\left( 
	%{{x_k} 
	%				+ 
	%						{s_k}} 
	%					\right)} \right)}}{{{m_k}\left( 0 \right) - {m_k}\left( {{s_k}} \right)}}} 
	%	\right|\le 
	%	\frac{{2{\varepsilon _h}\|s_k\|^2}}{{{m_k}\left( 0 \right) - {m_k}\left( {{s_k}} \right)}},$
	%\end{center}
	\begin{align*}
	\left| {{\rho _k} - {{\tilde \rho }_k}} \right| =\left| {\frac{{f\left( {{x_k}} \right) - 
				h\left( {{x_k}} \right) - \left( {f\left( {{x_k} + {s_k}} \right) - h\left( {{x_k} 
						+ 
						{s_k}} 
					\right)} \right)}}{{{m_k}\left( 0 \right) - {m_k}\left( {{s_k}} \right)}}} 
	\right|\le 
	\frac{{2{\varepsilon _h}\|s_k\|^2}}{{{m_k}\left( 0 \right) - {m_k}\left( {{s_k}} \right)}},
	\end{align*}
	that is
	\begin{align*}
	{{\tilde \rho }_k}- \frac{{2{\varepsilon _h}\|s_k\|^2}}{{{m_k}\left( 0 \right) - {m_k}\left( 
			{{s_k}} \right)}} \le 
	{\rho _k}  \le \frac{{2{\varepsilon _h}{\|s_k\|^2}}}{{{m_k}\left( 0 \right) - 
			{m_k}\left( {{s_k}} 
			\right)}}+ {{\tilde \rho }_k}.
	\end{align*}
	% $
	%{{\tilde \rho }_k}- \frac{{2{\varepsilon _h}\|s_k\|^2}}{{{m_k}\left( 0 \right) - {m_k}\left( 
	%		{{s_k}} \right)}} \le 
	%{\rho _k}  \le \frac{{2{\varepsilon _h}{\|s_k\|^2}}}{{{m_k}\left( 0 \right) - 
	%		{m_k}\left( {{s_k}} 
	%		\right)}}+ {{\tilde \rho }_k}.
	%$
	Then, define $\rho={{\tilde \rho }_k} - \frac{{2{\varepsilon 
				_h}{\|s_k\|^2}}}{{{m_k}\left( 0 \right) - 
			{m_k}\left( 
			{{s_k}} \right)}} $.  If $\rho \ge {\eta}$, we can obtain ${\rho _k} \ge {\eta}$. Thus, 
			in 
	the following analysis, we consider the size of $\rho$ that derive the desired lower bound 
	of the radius.

	Before giving the analyses, we briefly present the processing that why the $\Delta$ do not 
	approximate to zero. If  ${\rho} > {\eta 
	}$, the current iteration will be  successful and the radius $\Delta$ 
	will increase by a factor of $r_2$. Thus, we consider that whether there is a constant $C$ such 
	that
	$\Delta<C$ and the current iteration is successful simultaneously. Then we can see that such 
	constant $C$ is our desired bound of $\Delta$ due to the fact that $\Delta$ will increase again 
	under the successful iteration. What's more, such constant plays a critical role in determining 
	the iteration complexity.
	
	Instead of computing the $\rho$ directly, we consider another relationship, that is
	\begin{align}\label{Newton:STR:1-p}
	1 - \rho  = \frac{{{m_k}\left( 0 \right) - {m_k}\left( {{s_k}} \right) - \left( {h\left( 
				{{x_k}} \right) - h\left( {{x_k} + {s_k}} \right)} \right) + 2{\varepsilon 
				_h}{\|s_k\|^2}}}{{{m_k}\left( 0 \right) - {m_k}\left( {{s_k}} \right)}}
	\end{align}
	As long as $1-\rho<1-\eta$, we can see that $\rho>\eta$. Here, we consider the upper and lower 
	bound of denominator and numerator in (\ref{Newton:STR:1-p}) in Lemma 
	\ref{Newton:STR:Lemma:Upperbound_mx-msx} and Lemma \ref{Newton:STR:Lemma:Upperbound_fxs-msx}, 
	respectively. Moreover, we separately give the corresponding bound under the index sets 
	$\mathcal{S}_{\nabla f}$ and $\mathcal{S}_H$.
	%%%===================================================================
	%%%===================================================================
	%%%===================================================================
	\begin{lemma}\label{Newton:STR:Lemma:Upperbound_mx-msx}
		Suppose the Assumption 
		\ref{Newton:Assumption:approximation} and 
		\ref{Newton:Assumption:bound} hold,  $m_k(s)$ is defined in 
		(\ref{Newton:Definition_m}). For the case of
		$x_k\in\mathcal{S}_{\nabla f}$, if $m_k(s_k)\le m_k(s_k^C)$, where $s_k^C$ is the Cauchy 
		point, 
		then we 
		have
		\[{m_k}(0) - {m_k}({s_k}) \ge \frac{1}{2}({\varepsilon _{\nabla f}} - {\varepsilon _g})\min 
		\left\{ {{\Delta _k},({{{\varepsilon _{\nabla f}} - {\varepsilon _g}}})/{{{\kappa _H}}}} 
		\right\},\]
		For the case of $x_k\in \mathcal{S}_H$,  
		that is 
		${\lambda _{\text{min} }}( {\nabla^2f(x_k)} ) \le  - 
		{\varepsilon _H}$, 
		there exist a vector $s_k$ such that $\langle 
		{ g({x_k}),{s_k}} \rangle \le 0$, $ {s_k^T\nabla^2f(x_k){s_k}}<-v_{0}
		\Delta_k^2$, and $\left\| {{s_k}} \right\| = {\Delta 
			_k}$, where $v_{0}\ge\varepsilon_H$,	then we have
		\begin{align*}
		m_k( 0) - {m_k}(s_k )\ge \frac{1}{2} ( {{\varepsilon 
				_H} - {\varepsilon _B}})\Delta _k^2.
		\end{align*}
	\end{lemma}
	The solution $S_k$ for the \textit{Subproblem-Solver} is based on subproblem in 
	(\ref{Newton:Definition_m}). For the case of $x_k\in\mathcal{S}_{\nabla f}$, we use the Cauchy 
	point \cite{conn2000trust}; while for the case of $x_k\in\mathcal{S}_{H}$, there are many 
	methods 
	to derive the solution, such as Shift-and-invert \cite{garber2016faster},  Lanczos 
	\cite{kuczynski1992estimating} and Negative-Curvature \cite{carmon2018accelerated}.  We do not 
	present the 
	details information, which beyond our scope of this paper.
	%%%===================================================================
	%%%===================================================================
	%%%===================================================================
	\begin{lemma}\label{Newton:STR:Lemma:Upperbound_fxs-msx}
		Suppose  the Assumption 
		\ref{Newton:Assumption:function},
		\ref{Newton:Assumption:approximation} and \ref{Newton:Assumption:Bound:Variance} hold, 
		based 
		on the definition of $m_k(s)$ in 
		(\ref{Newton:Definition_m}),  we have
		\begin{align*}
		{m_k}\left( 0 \right) - {m_k}\left( {{s_k}} \right) - \left( {h\left( {{x_k}} \right) - 
			h\left( {{x_k} + {s_k}} \right)} \right) \le 2\left( 
		{\frac{{\mathbb{I}\left( 
					{\left| {{\mathcal{S}_h}} \right| < n} \right)}}{{\left| {{\mathcal{S}_h}} 
					\right|}}H_1 + 
			{\varepsilon _g}} 
		\right){\Delta _k} + \frac{3}{2}\left( {\frac{{\mathbb{I}\left( {\left| {{\mathcal{S}_h}} 
						\right| < 
						n} 
					\right)}}{{\left| {{\mathcal{S}_h}} \right|}}{H_2} + L_H{\Delta _k} + 
			{\varepsilon 
				_B}} 
		\right)\Delta 
		_k^2.
		\end{align*}
		If $\mathcal{S}_h=\mathcal{S}_g$, we have
		\begin{align*}
		{m_k}\left( 0 \right) - {m_k}\left( {{s_k}} \right) - \left( {h\left( {{x_k}} \right) - 
			h\left( {{x_k} + {s_k}} \right)} \right) \le\frac{3}{2}\left( {\frac{{\mathbb{I}\left( 
					{\left| {{\mathcal{S}_h}} \right| < 
						n} 
					\right)}}{{\left| {{\mathcal{S}_h}} \right|}}{H_2} + L_H{\Delta _k} + 
			{\varepsilon 
				_B}} 
		\right)\Delta 
		_k^2.
		\end{align*}
		%	where $s_k$ is obtained from the subproblem in 
		%(\ref{SFSNewton:STR:AlgorithmSTR-solution}).
	\end{lemma}
	
	Note that, in Algorithm 
	\ref{Newton:STR:Algorithm}, for the case of $\|g\left( {{x_k}} \right)\| \le {\varepsilon 	
		_{\nabla 	f}}+ {\varepsilon _g}$, that is 
	\begin{center}
		$\left\| {\nabla f\left( {{x_k}} \right)} \right\| \le 
		\left\| {g\left( {{x_k}} \right)} \right\| + \left\| {\nabla f\left( {{x_k}} \right) - 
		g\left( 
			{{x_k}} \right)} \right\| \le {\varepsilon _{\nabla f}} + {\varepsilon _g} + 
			{\varepsilon 
			_g} = 
		{\varepsilon _{\nabla f}} + 2{\varepsilon _g}$,
	\end{center}
	we set $\mathcal{S}_h=\mathcal{S}_g$. 
	$x\in\mathcal{S}_H$, in which $\left\| {\nabla f\left( {{x_k}} \right)} \right\| \le 
	{\varepsilon _{\nabla f}}$, satisfies such case\footnote{In the case of $\varepsilon_{\nabla 
			f}+2\varepsilon_{g}>\nabla f(x_k)>\varepsilon_{\nabla f}$, we have  
		$\mathcal{S}_h=\mathcal{S}_g$ 
		such that the equality (\ref{Newton:STR:Lemma:Bound:Radius-3}) in Appendix become 
		$\frac{{4{\varepsilon _g} + 
				\frac{3}{2}\left( {{L_H}{\Delta _k} + 2{\varepsilon _B} + 
					\frac{4}{3}{\varepsilon 
						_h}} \right){\Delta _k}}}{{\frac{1}{2}\left( {{\varepsilon _{\nabla f}} - 
					{\varepsilon _g}} 
				\right)}}$, which is smaller than equality (\ref{Newton:STR:Lemma:Bound:Radius-3}), 
		thus 
		$\Delta_{{\rm{min1}}}$ is also satisfying such case. In this paper, in order to simply the 
		analysis, we  consider the 
		case $\nabla f(x)>\varepsilon_{\nabla f}$ without the requirement of $\mathcal{S}_{\nabla 
			f}=\mathcal{S}_H$.}. 
	Thus, we give the $\Delta_{\text{min2}}$ in 
	(\ref{Newton:STR:Lemma:Bound:Radius-equality}) based on such implementation, which is key for 
	analyses. The reason we make such implementation is to 
	ensure that there is lower bound of radius. What's more, the parameters' setting is more 
	simple. 
	Based on above lemmas, we analyze the minimal radius.
	%===============================================================================
	%===============================================================================
	%===============================================================================
	\begin{lemma}\label{Newton:STR:Lemma:Bound:Radius}
		In Algorithm \ref{Newton:STR:Algorithm}, suppose the Assumption 
		\ref{Newton:Assumption:function}-
		\ref{Newton:Assumption:Bound:Variance} hold,  let $| {{\mathcal{S}_h}} | = \min 
		\{ n,\max \{ H_1/\varepsilon_g,H_2/\varepsilon_B 
		\} \}$, $1>r_1>0$, there 
		will be a  non-zero radium
		\begin{align}\label{Newton:STR:Lemma:Bound:Radius-equality}
		{\Delta _\text{min}} = \min \left\{ {{\Delta _\text{min1}},{\Delta 
				_\text{min2}}} 
		\right\},
		\end{align}		
		where the parameters satisfy
		\begin{align}
		\label{Newton:STR:Lemma:Bound:Radius-equality1}
		{\Delta _\text{min1}} =& {\kappa _1}\left( {{\varepsilon _{\nabla f}} - 
			{\varepsilon _g}} 
		\right),{\kappa _1} = {r_1}\min \left\{ {\frac{1}{{{\kappa _H}}},\frac{1}{{40}}\left( {1 - 
				\eta 
			} \right),\sqrt {\frac{1}{{12{L_H}}}\left( {1 - \eta } \right)} } \right\},\\
		\label{Newton:STR:Lemma:Bound:Radius-equality2}
		{\Delta _\text{min2}} =& {\kappa _2}\left( {{\varepsilon _H} - 
			{\varepsilon _B}} 
		\right),{\kappa _2} = {r_1}\frac{1}{{6{L_H}}}\left( {1 - \eta } \right),\\
		{\varepsilon _g} =& \frac{1}{{16}}\left( {1 - \eta } \right)\left( {{\varepsilon _{\nabla 
					f}} - {\varepsilon _g}} \right),\nonumber\\	
		{\varepsilon _B} =& {\varepsilon _h} = \frac{1}{{10}}\left( {1 - \eta } \right)\left( 
		{{\varepsilon _H} - {\varepsilon _B}} \right).\nonumber
		\end{align}
		In particular, 	$\Delta_\text{min1}$ belongs to the case of $x\in 
		\mathcal{S}_{\nabla f}$ and $\Delta_\text{min2}$ belongs to the case of 
		$x\in 
		\mathcal{S}_{H}$.
	\end{lemma}
	%%%%%%%%%%%%%%%%%%%%%%%%%%%%%%%%%%%%%%%%%%%%%%%%%%%%%%%%%%%%%%%%%%%%%%%%%%%%%%%%%%%%%%%
	%%%%%%%%%%%%%%%%%%%%%%%%%%%%%%%%%%%%%%%%%%%%%%%%%%%%%%%%%%%%%%%%%%%%%%%%%%%%%%%%%%%%%%%
	%%%%%%%%%%%%%%%%%%%%%%%%%%%%%%%%%%%%%%%%%%%%%%%%%%%%%%%%%%%%%%%%%%%%%%%%%%%%%%%%%%%%%%%
	\subsection{Convergence and iteration complexity}
	In this section, we  present the successful and 
	unsuccessful iteration complexity based on Lemma 
	\ref{Newton:STR:Lemma:Bound:Radius} including 
	$x\in \mathcal{S}_{\nabla f}$ and $k\in \mathcal{S}_H$, and then provide the total number of 
	iterations.
	%%%=============================================================
	%%%=============================================================
	%%%=============================================================
	\begin{theorem}\label{Newton:STR:Theorem:Iteration}
		In Algorithm \ref{Newton:STR:Algorithm}, suppose the Assumption 
		\ref{Newton:Assumption:function}-
		\ref{Newton:Assumption:Bound:Variance} hold, let $| {{\mathcal{S}_h}} | = \min 
		\{ n,\max \{ H_1/\varepsilon_g,H_2/\varepsilon_B 
		\} \}$, $\{ {f( {{x_k}} )} \}$ is 
		bounded below by $f_{\text{low}}$, the number of successful iterations ${T_\text{suc}}$  is 
		no large than
		\[{\kappa _3}{\rm{max}}\{ {{{( {{\varepsilon _{\nabla f}} - {\varepsilon _g}} 
					)}^{ - 2}},{{({\varepsilon _H} - {\varepsilon _B})}^{ - 3}}} \},\]
		where ${\kappa _3} = 2\left( {f\left( {{x_0}} \right) - {f_{low}}} \right){\rm{max}}\left\{ 
		{1/\left( {\eta {\kappa _1}} \right),1/\left( {\eta \kappa _2^2} \right)} \right\}$, 
		$\kappa_1$ 
		and $\kappa_2$ are defined in (\ref{Newton:STR:Lemma:Bound:Radius-equality1}) and 
		(\ref{Newton:STR:Lemma:Bound:Radius-equality2}). The number of unsuccessful iterations 
		$T_\text{unsuc}$ is at 
		most 
		\[ \frac{1}{{ - \log {r_1}}}\left( {\log \left( {\frac{{{\Delta 
							_{{\text{max}}}}}}{{{\Delta _{{\text{min}}}}}}} \right) - T\log {r_2}} 
		\right).\]
		where ${\Delta _{\text{max}}}$ and ${\Delta _{\text{min} }}$ are defined in  
		(\ref{Newton:STR:Lemma:Bound:Radius-equality}) and Algorithm \ref{Newton:STR:Algorithm}, 
		$1>{r}_1>0$ and 
		$r_2\ge 
		1$. Thus, the total 
		number of iterations is 
		\[\mathcal{O}\left( {\max \{ {{{({\varepsilon _{\nabla f}} - {\varepsilon _g})}^{ - 
						2}},{{({\varepsilon _H} - {\varepsilon _B})}^{ - 3}}} \}} \right)\]
		%	\begin{align*}
		%	\max \{ {\mathcal{O}( {{{({\varepsilon _{\nabla f}} - {\varepsilon _g})}^{ 
		%					- 2}}{{({\varepsilon _H} - {\varepsilon _B})}^{ - 1}}} ),\mathcal{O}( 
		%		{{{({\varepsilon _H} - {\varepsilon _B})}^{ - 3}}} )} \}.
		%	\end{align*}
	\end{theorem}

	After the $|T_\text{suc}|+|T_\text{unsuc}|$ iterations, it will fall into $\mathcal{S}_*$ and 
	converge to the 
	stationary point. As can be seen from above Theorems, we make two conclusions: the first is the 
	order of iteration complexities are the same as \cite{cartis2011adaptiveb} and 
	\cite{xu2017second} 
	if the parameters $\varepsilon_g$ and $\varepsilon_B$ are set properly according to 
	$\varepsilon_{\nabla f}$ and $\varepsilon_H$ respectively; 
	the second is that when $| {{\mathcal{S}_h}} |<n$, the total number of computation iteration 
	including computing the function is less than that of \cite{cartis2011adaptiveb} and 
	\cite{xu2017second}; when  $| {{\mathcal{S}_h}} |=n$, our result is equal 
	to\cite{cartis2011adaptiveb} and 
	\cite{xu2017second}. Thus, our proposed algorithm is more general. 
	%%%%%%%%%%%%%%%%%%%%%%%%%%%%%%%%%%%%%%%%%%%%%%%%%%%%%%%%%%%%%%%%%%%%%%%%%%%%%
	%%%%%%%%%%%%%%%%%%%%%%%%%%%%%%%%%%%%%%%%%%%%%%%%%%%%%%%%%%%%%%%%%%%%%%%%%%%%%
	%%%%%%%%%%%%%%%%%%%%%%%%%%%%%%%%%%%%%%%%%%%%%%%%%%%%%%%%%%%%%%%%%%%%%%%%%%%%%
	
	\section{Stochastic Adaptive Regularization using Cubics}\label{Newton:Section:SARC}
	
	\begin{algorithm*}[t]
		\caption{SARC  with the Inexact Hessian and Gradient}
		\label{Newton:SARC:Algorithm}
		\begin{algorithmic}[1]
			\Require given $x_0$, $r_2\ge 1> r_1$, $1>\eta>0$, $\varepsilon>0$ and  
			$\sigma_\text{min}>0$.
			\Ensure
			\For{$k$=1 to T }
			\State Compute the approximate gradient $g(x_k)$ and 
			Hessian matrix 	$B(x_k)$ based on 
			(\ref{Newton:STR:Subsampling:g}) 
			and (\ref{Newton:STR:Subsampling:B}).
			\State \Comment({\r1 if 
				%${\lambda _{\text{min} }}( {\nabla^2f(x_k)} ) 	\le  - 	{\varepsilon _H}$, 
				%$||\nabla f\left( {{x_k}} \right)|| \le {\varepsilon _{\nabla f}}$, 
				$||g\left( {{x_k}} \right)|| \le {\varepsilon 	_{\nabla 	f}} + {\varepsilon 
				_g}$, 
				we set $\mathcal{S}_h=\mathcal{S}_g$})
			\State Compute the direction vector $s_k$
			\begin{align}\label{Newton:SARC:Algorithm-solution}
			s_k=\text{Subproblem-Solver}(g(x_k),B(x_k),\sigma_k).
			\end{align}	
			\State Compute the approximate function $h(x_k)$ and $h(x_k+s_k)$ based on 
			(\ref{Newton:STR:Subsampling:h}).	
			\State 	Compute $
			{{\tilde \rho }_k} = \frac{{h\left( {{x_k}} \right) - h\left( {{x_k} + {s_k}} 
					\right)}}{{{p_k}\left( 0 \right) - {p_k}\left( {{s_k}} \right)}}
			$
			\State Set $\rho={{{\tilde \rho }_k} - \frac{{2{\varepsilon _f}/{\r1\sigma^ 
							2_k}}}{{{p_k}\left( 0 
						\right) - {p_k}\left( {{s_k}} \right)}}}$			
			\State set ${x_{k + 1}} = \left\{ {\begin{array}{*{20}{l}}
				{{x_k} + {s_k},}&{\rho  \ge \eta ,}\\
				{{x_k},}&{{\rm{otherwise}}{\rm{.}}}
				\end{array}} \right.$
			\State set
			${\sigma _{k + 1}} = \left\{ {\begin{array}{*{20}{l}}
				{\max \left\{ {{\sigma _{\text{min} }},{r_1}{\sigma _k}} \right\},}&{\rho  > \eta 
				,}\\
				{{r_2}{\sigma _k},}&{{\rm{otherwise}}.}
				\end{array}} \right.$	
			\EndFor
		\end{algorithmic}
	\end{algorithm*}
	
	In this section, we consider the stochastic adaptive regularization using 
	Cubics (SARC) method, which solves the following unconstrained minimization 
	problem at each iteration: 
	\begin{align}\label{Newton:SARC:definition_P}
	\mathop {\min }\nolimits_{s\in \mathbb{R}^d}  p_k( s ) := m_k(s) + \frac{{{\sigma _k}}}{3}{\| s 
		\|^3},
	\end{align}
	where $m_k(s)$ is defined in (\ref{Newton:Definition_m}), and $\sigma_k$ 
	is an adaptive parameter that can be considered as the reciprocal of the 
	trust-region radius. Algorithm \ref{Newton:SARC:Algorithm} presents the process for 
	updating the $x_k$ and $\sigma_k$. Similar to the analysis as in \cite{cartis2011adaptivea}, 
	$\sigma_k$ in the 
	cubic term actually performs one more task, besides accounting for the 
	discrepancy between the objective function and its corresponding 
	second-order Taylor expansion, but also for the difference between   the exact and approximate 
	function, gradient and Hessian. The update rules of $\sigma$ is analogous to stochastic region 
	method. $\sigma$ will decrease if sufficient decrease is obtained in some measure of relative 
	objective chance, but increase otherwise. Following the framework of STR, we analyze SARC 
	including 
	two parts: To ensure the existence of the maximal bound of $\sigma$ and present the iterative 
	complexity.
	
	%%%%%%%%%%%%%%%%%%%%%%%%%%%%%%%%%%%%%%%%%%%%%%%%%%%%%%%%%%%%%%%%%%%%%%%%%%%%%%%%%%%%%%%
	%%%%%%%%%%%%%%%%%%%%%%%%%%%%%%%%%%%%%%%%%%%%%%%%%%%%%%%%%%%%%%%%%%%%%%%%%%%%%%%%%%%%%%%
	%%%%%%%%%%%%%%%%%%%%%%%%%%%%%%%%%%%%%%%%%%%%%%%%%%%%%%%%%%%%%%%%%%%%%%%%%%%%%%%%%%%%%%%
	\subsection{Bounds analysis of  the adaptive parameter}
	Similar to STR, we present the definition of $\rho$ directly,
	\begin{align}\label{Newton:SARC:1-p}
	1 - \rho  =& \frac{{{p_k}\left( 0 \right) - {p_k}\left( {{s_k}} \right) - \left( {h\left( 
				{{x_k}} \right) - h\left( {{x_k} + {s_k}} \right)} \right) + 2{\varepsilon 
				_h}\|s_k\|^2}} {{{p_k}\left( 0 \right) - {p_k}\left( {{s_k}} \right)}}
	\end{align}
	In order to satisfy inequality (\ref{Newton:SARC:1-p}), we 
	need to obtain the lower bound of the numerator in 
	(\ref{Newton:SARC:1-p}). Firstly, we derived the lower bound 
	from the view of a Cauchy pint, but subject to ${p_k}\left( {{s_k}} \right) 
	\le {p_k}\left( {s_k^C} \right)$. Note that the lower bound is almost the 
	same 
	as in \cite{cartis2011adaptivea} and \cite{xu2017newton}, but give the 
	proof 
	from the geometrical explanation.
	%%%===========================================
	%%%===========================================
	%%%===========================================
	\begin{lemma}\label{Newton:SARC:Lemma:UpperboundOfP0Ps} 
		Suppose that the step size $s_k$ satisfying ${p_k}\left( {{s_k}} 
		\right) 
		\le {p_k}\left( {s_k^C} \right)$, where $s_k^C$ is a Cauchy point, 	defined 
		as
		\begin{align*}
		s_k^C =&  - {\alpha _k}g\left( {{x_k}} \right),
		%\\
		{\alpha _k} =
		%& 
		\mathop {\arg \min }\nolimits_{\alpha  \in {\mathbb{R}_ + 
		}} 
		\left\{ {{p_k}\left( x_k{ - \alpha g\left( {{x_k}} \right)} \right)} 
		\right\},
		\end{align*}
		for all $k\ge 0$,
		we have that 
		%      $p_k(0) - {p_k}( s_k)
		% \ge 
		% \frac{1}{10}\left\| {{g_k}} \right\|\min \left\{ 
		% {{\left\| {{g_k}} \right\|}/{{\left\| {{B_k}} \right\|}},\sqrt {{{\left\| {{g_k}} 
		% 				\right\|}}/(\sigma _k)} } 
		% \right\}.$
		\begin{center}
			$p_k(0) - {p_k}( s_k)
			\ge 
			\frac{1}{10}\left\| {{g_k}} \right\|\min \left\{ 
			{{\left\| {{g_k}} \right\|}/{{\left\| {{B_k}} \right\|}},\sqrt {{{\left\| {{g_k}} 
							\right\|}}/(\sigma _k)} } 
			\right\}.$
		\end{center}
		Specifically, we set ${\alpha} = 2/(\| B_k
		\| + \sqrt 
		{{\| B_k \|}^2 + 4{\sigma _k}\| g_k 
			\|} )$ and $s_k=-\alpha g(x_k)$ that satisfy above inequality. Furthermore, we 
		can also obtain the upper bound of the step 
		$\| s_k \|$, $k>0$, satisfies 
		$\| s_k \| \le {{11}}/{4}\max \{ {{{\| 
					{B({x_k})} \|}}/{{{\sigma _k}}},\sqrt {{{\| 
						{g({x_k})} 
						\|}}/{{{\sigma _k}}}} } \}.$		
		%		\begin{center}
		%			$\| s_k \| \le {{11}}/{4}\max \{ {{{\| 
		%						{B({x_k})} \|}}/{{{\sigma _k}}},\sqrt {{{\| 
		%							{g({x_k})} 
		%							\|}}/{{{\sigma _k}}}} } \}.$
		%		\end{center}
	\end{lemma}
	
	Following the subspace analysis in the cubic model as in 
	\cite{cartis2011adaptivea} 
	and \cite{cartis2011adaptiveb}, in order to widen the scope of convergence 
	analysis and iteration complexity, we also consider the step size $s_k$ on 
	the 
	following conditions
	\begin{align}
	\label{Newton:SARC:Assumption_s1}
	&\left\langle {g\left( {{x_k}} \right),{s_k}} \right\rangle  + 
	s_k^T{B_k}{s_k} + 
	{\sigma _k}{\left\| {{s_k}} \right\|^3} = 0,\\
	\label{Newton:SARC:Assumption_s2}
	&s_k^T{B_k}{s_k} + {\sigma _k}{\left\| {{s_k}} \right\|^3} \ge 0,\\
	\label{Newton:SARC:Assumption_s3}
	&\left\| {\nabla {p_k}\left( {{s_k}} \right)} \right\| \le {\theta _k}\left\| {\nabla g\left( 
		{{x_k}} \right)} \right\|,{\theta _k} \le {\kappa _\theta }\text{min}\left\{ {1,\left\| 
		{{s_k}} 
		\right\|} \right\},{\kappa _\theta } < 1.
	\end{align}
	Thus, we can also obtain lower bound of the numerator in 
	(\ref{Newton:SARC:1-p}), which will be used for analyzing the 
	convergence and 
	iteration complexity in the case of the saddle point.
	%========================================================================
	%=========================================================================
	%=====================================================================
	\begin{lemma}\label{Newton:SARC:Lemma:UpperboundOfP0Ps-2} Given the  
		conditions of 
		$s_k$ in 
		(\ref{Newton:SARC:Assumption_s1}) and 
		(\ref{Newton:SARC:Assumption_s2}), 
		we have
		%$p( 0) - {p_k}( s_k ) \ge \frac{{{\sigma_k}}}{6}{\| {{s_k}} \|^3}.$
		\begin{align*}
		p( 0) - {p_k}( s_k ) \ge 
		\frac{{{\sigma_k}}}{6}{\| {{s_k}} \|^3}.
		\end{align*}
		Furthermore, 	suppose 
		Assumption \ref{Newton:Assumption:function} and	
		\ref{Newton:Assumption:approximation},  and the condition (\ref{Newton:SARC:Assumption_s3}) 
		hold, for $\left\| {g\left( {{x_{k + 1}}} \right)} \right\| \ge {\varepsilon _{\nabla f}} - 
		{\varepsilon _g}$, we have
		\begin{align}\label{Newton:SARC:Lemma:UpperboundOfP0Ps-2-1}
		\left\| {g\left( {{x_{k + 1}}} \right)} \right\| \le {\kappa 
			_s}{\left\| {{s_k}} 
			\right\|^2},
		\end{align}
		where 
		\[\begin{array}{l}
		{\kappa _s} = \min \left\{ {\frac{{2{\varepsilon _B} + \left( 
					{{L_H} + {\sigma _k}} 
					\right) + 2{\kappa _\theta }{\varepsilon _g} + {\kappa _\theta 
					}{L_{\nabla f}}}}{{\left( {1 - 
						{\theta 
							_k}} \right)}},\frac{{{L_H} + {\sigma _k} + {\kappa _\theta }{L_{\nabla 
							f}}}}{{1 - {\theta _k} 
					- {\zeta 
						_1} - {\zeta _2}}}} \right\}, {\zeta _1},{\zeta _2} < 1\\
		{\varepsilon _B} \le {\zeta _1}\left( {{\varepsilon _{\nabla f}} - {\varepsilon _g}} 
		\right),{\varepsilon _g} \le {\zeta _2}\left( {{\varepsilon _{\nabla f}} - {\varepsilon 
		_g}} 
		\right)
		\end{array}\]
	\end{lemma}
	Different from the STR, we can also derive the relationship between $g(x_{k+1})$ and $\|s_k\|$. 
	The 
	core process of the proof is based on cubic regularization of Newton method 
	\cite{nesterov2006cubic}. Such a relationship leads to the improved iteration complexity. 
	Besides 
	lower bound of the numerator in (\ref{Newton:SARC:1-p}), we can also obtain the corresponding 
	upper bound of the denominator, which is similar to Lemma 
	\ref{Newton:STR:Lemma:Upperbound_fxs-msx}. 
	Thus, we do not provide the proof.
	%%%====================================================================
	%%%====================================================================
	%%%====================================================================
	\begin{lemma}\label{Newton:SARC:Lemma:Upperbound_p-h}
		Suppose the Assumption \ref{Newton:Assumption:function}, 
		\ref{Newton:Assumption:approximation} and \ref{Newton:Assumption:Bound:Variance} hold, at 
		$k$-iteration, we have
		\begin{align*}
		&{p_k}\left( 0 \right) - {p_k}\left( {{s_k}} \right) - \left( {h\left( {{x_k}} \right) - 
			h\left( 
			{{x_k} + {s_k}} \right)} \right)\\
		\le& 2\left( {\frac{\mathbb{I}{\left( {\left| {{{\cal S}_h}} \right| < n} \right)}}{{\left| 
					{{{\cal 
								S}_h}} 
					\right|}}{H_1} + {\varepsilon _g}} \right)\left\| {{s_k}} \right\| + 
		\frac{3}{2}\left( 
		{\frac{\mathbb{I}{\left( {\left| {{{\cal S}_h}} \right| < n} \right)}}{{\left| {{{\cal 
								S}_h}} 
					\right|}}{H_2} + {\varepsilon _B}} \right){\left\| {{s_k}} \right\|^2} + \left( 
		{\frac{3}{2}{L_H} - \frac{1}{3}{\sigma _k}} \right){\left\| {{s_k}} \right\|^3}
		\end{align*}
		If $\mathcal{S}_h=\mathcal{S}_g$, we have
		\begin{align*}
		{p_k}\left( 0 \right) - {p_k}\left( {{s_k}} \right) - \left( {h\left( {{x_k}} \right) - 
			h\left( 
			{{x_k} + {s_k}} \right)} \right)\le\frac{3}{2}\left( 
		{\frac{\mathbb{I}{\left( {\left| {{{\cal S}_h}} \right| < n} \right)}}{{\left| {{{\cal 
		S}_h}} 
					\right|}}{H_2} + {\varepsilon _B}} \right){\left\| {{s_k}} \right\|^2} + \left( 
		{\frac{3}{2}{L_H} - \frac{1}{3}{\sigma _k}} \right){\left\| {{s_k}} \right\|^3}.
		\end{align*}
	\end{lemma}
	%%========================================
	%%===========================================
	Based on the above lemmas, we can derive the upper bound of adaptive parameter $\sigma$, which 
	is 
	used to analyze the iteration complexity. Furthermore, the parameters' setting, such as 
	$\varepsilon_g$, $\varepsilon_B$ and $\varepsilon_h$ are similar 
	to that of Lemma \ref{Newton:STR:Lemma:Bound:Radius}.
	\begin{lemma}\label{Newton:SARC:Lemma:bound_sigmal}In Algorithm \ref{Newton:SARC:Algorithm}, 
		suppose Assumption  
		\ref{Newton:Assumption:function}-\ref{Newton:Assumption:Bound:Variance} hold, let $| 
		{{\mathcal{S}_h}} | = \min 
		\{ n,\max \{ H_1/\varepsilon_g,H_2/\varepsilon_B 
		\} \}$, $r_2>1$, the parameter 
		$\sigma$ is bounded by 
		\begin{align}\label{Newton:SARC:Lemma:bound_sigmal-1}
		{\sigma _{\text{max}}} = \max  \left\{ {{\sigma _\text{max1}},{\sigma 
				_\text{max2}}} 
		\right\},
		\end{align}	
		where
		\begin{align}
		{\sigma _{{\text{max1}}}} =& {\kappa 
			_4}\frac{1}{{({\varepsilon _{\nabla f}} - {\varepsilon _g})}},{\sigma _{\text{max2}}} = 
		\frac{9}{2}{r_2}{L_H},\nonumber\\
		\label{Newton:SARC:Lemma:bound_sigmal-2}
		{\kappa _4} =& {r_2}\left\{ {\kappa _H^2,\frac{{{{\left( {304\left( {3{\varepsilon _B} + 
									2{\varepsilon _h}} \right)} \right)}^2}}}{{\left( {1 - \eta } 
					\right)}},\frac{9}{2}({\varepsilon _{\nabla f}} - {\varepsilon _g}){L_H}} 
		\right\},\\
		{\varepsilon _g} =& 
		\frac{1}{{220}}\left( {1 - \eta } \right)({\varepsilon _{\nabla f}} - 
		{\varepsilon _g}),\nonumber\\				
		{\varepsilon _B} =& {\varepsilon _h} =\frac{1}{{36}}\left( {1 - \eta } \right)\left( 
		{{\varepsilon _H} - 
			{\varepsilon _B}} \right).\nonumber
		\end{align}
		In particular, 	$\sigma_\text{min1}$ belongs to the case of $x\in 
		\mathcal{S}_{\nabla f}$ and $\sigma_\text{min2}$ belongs to the case of 
		$x\in 
		\mathcal{S}_{H}$.	
	\end{lemma}

	%%%%%%%%%%%%%%%%%%%%%%%%%%%%%%%%%%%%%%%%%%%%%%%%%%%%%%%%%%%%%%%%%%%%%%%%%%%%%%%%%%%%%%%
	%%%%%%%%%%%%%%%%%%%%%%%%%%%%%%%%%%%%%%%%%%%%%%%%%%%%%%%%%%%%%%%%%%%%%%%%%%%%%%%%%%%%%%%
	%%%%%%%%%%%%%%%%%%%%%%%%%%%%%%%%%%%%%%%%%%%%%%%%%%%%%%%%%%%%%%%%%%%%%%%%%%%%%%%%%%%%%%%
	\subsection{Analysis of convergence and iteration complexity}
	Based on above lemmas, we present the iteration complexity. Different from STR, we derive two 
	kinds of complexity. The first one has the same order as STR while the second is better or 
	equal to 
	STR. The difference lies that if the criterion conditions in (\ref{Newton:SARC:Assumption_s3}) 
	is 
	satisfied, the iteration 
	complexity to the stationary point is improved.
	%%%==============================================================
	%%%============================================================
	%%%================================================================
	\begin{theorem}\label{Newton:SARC:theorem:Iteration}
		In Algorithm \ref{Newton:SARC:Algorithm}, suppose the Assumption 
		\ref{Newton:Assumption:function}-
		\ref{Newton:Assumption:Bound:Variance} hold, let $| {{\mathcal{S}_h}} | = \min 
		\{ n,\max \{ H_1/\varepsilon_g,H_2/\varepsilon_B 
		\} \}$, $\{ {f( {{x_k}} )} \}$ is 
		bounded below by $f_{\text{low}}$, the number of successful iterations ${T_\text{suc}}$  is 
		no large than
		\[{\kappa _5}\max \{ {{{({\varepsilon _{\nabla f}} - {\varepsilon _g})}^{ - 2}},{{( 
					{{\varepsilon _H} - {\varepsilon _B}} )}^{ - 3}}} \},\]
		where ${\kappa _5} = ( {f( {{x_0}} ) - {f_{{\rm{low}}}}} )\max \{ 
		{5/( \eta \kappa _4^{ - 1/2} ),6{\sigma _{{\rm{max2}}}}/\eta } \}$, 
		$\kappa_4$  is defined in (\ref{Newton:SARC:Lemma:bound_sigmal-2}). If the conditions 
		(\ref{Newton:SARC:Assumption_s1})-(\ref{Newton:SARC:Assumption_s3}) are satisfied, then the 
		number of successful 	iterations 
		$T_\text{unsuc}$ is at 
		most 
		\[{\kappa _6}\max \{ {{{({\varepsilon _{\nabla f}} - {\varepsilon _g})}^{ - 3/2}},{{( 
					{{\varepsilon _H} - {\varepsilon _B}} )}^{ - 3}}} \},\]
		where ${\kappa _6} = ( {f( {{x_0}} ) - {f_{{\text{low}}}}} )\max \{ 
		{6\kappa _s^{3/2}/( {\eta {\sigma _{{\text{min}}}}} ),6{\sigma _{{\text{max2}}}}/\eta } 
		\}$, $\kappa_s$ is defined (\ref{Newton:SARC:Lemma:UpperboundOfP0Ps-2-1}).
	\end{theorem}
	As can be seen above results, our proposed method has the same order of iteration complexity as 
	in 
	\cite{xu2017newton} and \cite{xu2017second}. However, our algorithm does not require the full 
	computation of function and gradient with the finite-sum structure such that reduce the 
	computation 
	cost properly.
	%%%%%%%%%%%%%%%%%%%%%%%%%%%%%%%%%%%%%%%%%%%%%%%%%%%%%%%%%%%%%
	%%%%%%%%%%%%%%%%%%%%%%%%%%%%%%%%%%%%%%%%%%%%%%%%%%%%%%%%%%%%%%
	%%%%%%%%%%%%%%%%%%%%%%%%%%%%%%%%%%%%%%%%%%%%%%%%%%%%%%%%%%%%
	%%%%%%%%%%%%%%%%%%%%%%%%%%%%%%%%%%%%%%%%%%%%%%%%%%%%%%%%%%%
	
	\section{Experiment}\label{Newton:Section:Experiment}
	
	\begin{figure*}[t]
		\centering
		\includegraphics[width=0.7\linewidth]{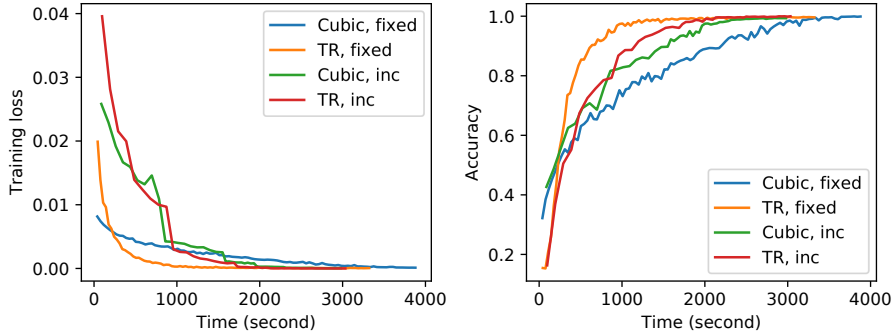}
		\caption{Training loss/accuracy vs running time for different schemes of 
			second-order methods. Note that TR-full and ARC-full only finish 1 epoch 
			after 3000 seconds, so we omit them in the plot. }
		\label{fig:experiment_result}
	\end{figure*}

%	\begin{figure*}[t]
%	\centering
%	\includegraphics[width=0.42\textwidth]{indexset2}
%	\caption{Illustration of three situations in analyzing the convergence 
%		to first and second critical 
%		point.}\label{Newton:Figure3Situation}
%\end{figure*}

	In this section, we give a comprehensive comparison among trust region and ARC 
	algorithms. Our goal is not to show TR/ARC methods are state-of-the-art 
	compared with other solvers; we are trying to present different variances of TR and 
	ARC, and show using a fixed batch size to estimate both gradient and Hessian
	is the best choice for large-scale optimization. 
	We compare following variants:
	\begin{itemize}
		\item Stochastic TR with fixed batch size (TR, fixed): Both Hessian and 
		gradient are estimated through a fixed batch size at each iteration. In our 
		experiment, we choose batch size $|B_g|=|B_H|=256$. 
		\item Stochastic TR with growing batch size (TR, inc): As above, both 
		Hessian and gradient are estimated through a batch of samples, except 
		that the sample size is increasing with epochs: At the early stage we feed 
		a crude estimation of $\nabla f(x)$ and $\nabla^2 f(x)$ on small batch, 
		then gradually increase the batch size to give better gradient and Hessian 
		information. In practice we multiply the batch size by a factor of $2$ for 
		every $10$ epochs until memory is used up. 
		\item Stochastic TR with exact gradient and subsampled Hessian (TR, full): We use 
		implementation similar with \cite{xu2017newton,xu2017second}. At each iteration the 
		gradient is 
		exact, while the Hessian is approximated on a batch of $256$ samples.
		\item Stochastic ARC with fixed batch size (Cubic, fixed): This is similar to stochastic 
		TR, 
		the batch size is fixed to $|B_g|=|B_H|=256$.
		\item Stochastic ARC with growing batch size (Cubic, inc): Similar to stochastic TR with 
		increasing batch size, we double the batch size once every $10$ epochs to estimate both 
		Hessian 
		and gradient.
		\item ARC with exact gradient and subsampled Hessian (Cubic, full): The 
		gradient is exactly computed while the Hessian is approximated on a batch 
		of $256$ samples. 
	\end{itemize}
	Note that the ARC-inc algorithm is proposed in~\cite{kohler2017sub}; TR-inc is 
	a generalization of that to the trust region case; TR-full and ARC-full are 
	proposed in~\cite{xu2017newton}; TR-fixed and ARC-fixed are our method 
	analyzed in this paper. 
	
	We train a VGG16 network\footnote{Publically available at 
		\url{https://raw.githubusercontent.com/kuangliu/pytorch-cifar/master/models/vgg.py}} on 
		CIFAR10 
	dataset using the above algorithms and evaluate the performance according to
	training loss and test accuracy with respect to training time. All the 
	experiments are run on a machine with 1 Titan Xp GPU. 
	The result is reported in Figure~\ref{fig:experiment_result}.
	
	Note that TR-full and ARC-full take 3500 seconds  while 
	TR-fixed and ARC-fixed only take 40 seconds for each epoch. Thus we omit the 
	full gradient versions on the plot since there will be only one point there.
	This shows that calculating the full gradient for each update is 
	too expensive for solving large-scale problems like deep learning. 
	Apart from that we notice both Cubic and TR on fixed batch size are faster 
	than their growing batch size version, this validates our guess that sampling 
	a fixed number of data to estimate gradient and Hessian is sufficient to make 
	trust region and ARC work, and this scheme turns out to be more efficient than 
	sample a growing batch over time. Moreover, in our task, the trust region 
	algorithm is faster than ARC algorithm. However, we are not sure whether this 
	phenomenon also applies to other tasks.

	%%%%%%%%%%%%%%%%%%%%%%%%%%%%%%%%%%%%%%%%%%%%%%%%%%%%%%%%%%%%%%%%%%%%%%%%%%%%%
	%%%%%%%%%%%%%%%%%%%%%%%%%%%%%%%%%%%%%%%%%%%%%%%%%%%%%%%%%%%%%%%%%%%%%%%%%%%%%
	%%%%%%%%%%%%%%%%%%%%%%%%%%%%%%%%%%%%%%%%%%%%%%%%%%%%%%%%%%%%%%%%%%%%%%%%%%%%%
	
	%\input{Section-STR}
	%\input{Section-SARC}
	%\input{Section-Exp}
	\section{Conclusion}\label{Newton:Section:Conclusion}
	In this paper, we present a family of stochastic trust region method and 
	stochastic cubic regularization method under inexact gradient and Hessian 
	matrix. Furthermore, in order to reduce the computation cost for the function 
	value of $f(x)$ in evaluating the role of $\rho$, we also present a 
	sub-sample technique to estimate the function.  We provide the theoretical 
	analysis of convergence and iteration complexity and obtain that we keep the 
	same order of iteration complexity but reduce the computation cost per 
	iteration. We apply our proposed method to deep learning application which 
	outperforms  the previous second-order methods. 
	{\small
		\bibliographystyle{unsrt}
		\bibliography{STR_bib}
	}
	
	%===========================================
	%==============\Appendix=====================
	\appendix

	\section{Proof of Sub-sampling}
	\textbf{Proof of Lemma \ref{Newton:Appendix:Tool:RandomSubset}}
	\begin{proof} Based on the  $\sum\nolimits_{i = 1}^m {{v_i}}  = \vec 0$, and 
		permutation and combination, 
		%\begin{itemize}
		%	\item 
		For the case that $\cal A$ is a non-empty, 
		uniformly random subset of $[m]$, we 
		have 
		\begin{align*}
		%&
		{\mathbb{E}_{\cal A}} {{{\left\| {\sum\nolimits_{b \in \mathcal{A}} 
						{{v_b}} } 	\right\|}^2}} 
		%\\ 
		=
		& 
		{\mathbb{E}_{\cal 
				A}}\left[ {\sum\nolimits_{b 
				\in {\cal A}} 	{{{\left\| 	{{v_b}} \right\|}^2}} } \right] 
		+ 
		\frac{1}{{C_n^A}}\sum\limits_{i 	\in [n]} 
		{\left\langle {{v_i},\frac{{C_{n - 1}^{A - 1}\left( {A - 1} 
						\right)}}{{n - 
						1}}\sum\limits_{i \ne j} {{v_j}} } \right\rangle 
		} 
		\\
		%\label{SCSG-Zero:Appendix:lemma-expectationSubset-equality1}
		=& A\frac{1}{n}\sum\nolimits_{i = 1}^n {v_i^2}  + \frac{{A\left( {A - 
					1} 
				\right)}}{{n\left( {n - 1} \right)}}\sum\nolimits_{i \in 
			[n]} 
		{\left\langle 
			{{v_i},\sum\nolimits_{i \ne j} {{v_j}} } \right\rangle } \\
		=& A\frac{1}{n}\sum\nolimits_{i = 1}^n {v_i^2}  + \frac{{A\left( {A - 
					1} 
				\right)}}{{n\left( {n - 1} \right)}}\sum\nolimits_{i \in 
			[n]} 
		{\left\langle 
			{{v_i}, - {v_i}} \right\rangle }\\ 
		%\label{SCSG-Zero:Appendix:lemma-expectationSubset-equality2}
		=& \frac{{A\left( {n - A} \right)}}{{\left( {n - 1} 
				\right)}}\frac{1}{n}\sum\nolimits_{i = 1}^n {v_i^2}
		\le A\mathbb{I}\left( {A < n} \right)\frac{1}{n}\sum\nolimits_{i = 
			1}^n 
		{v_i^2}. 
		\end{align*}
		Thus, we have
		%	\begin{align*}
		${\mathbb{E}_{\cal A}} {{{\left\| \frac{1}{A}{\sum\nolimits_{b \in 
							\mathcal{A}} 
						{{v_b}} } 	\right\|}^2}}=\frac{1}{A^2}{\mathbb{E}_{\cal 
				A}} {{{\left\| 
					{\sum\nolimits_{b \in \mathcal{A}} 
						{{v_b}} } 	\right\|}^2}}\le\frac{{\mathbb{I}\left( {A < 
					n} 
				\right)}}{A}\frac{1}{n}\sum\limits_{i = 1}^n {v_i^2}.
		$
		%	\end{align*}	
	\end{proof}	
	%%%=============================================================================================
	%%%=============================================================================================
	%%%=============================================================================================
	\textbf{Proof of Lemma \ref{Newton:Subsample:Lemma} }
	\begin{proof}
		Let us define $
		{X_i} = {f_i}( x ) - f( x ),i \in {\mathcal{S}_h}.$
		Based on the Assumption 
		\ref{Newton:Assumption:bound}, we have
		\begin{align*}
		| {{X_i}} | = | {{f_i}( x ) - f( x 
			)} | \le | {{f_i}( x )} | + 
		| {f( x )} | \le 2{\kappa _f} \Rightarrow 
		{| {{X_i}} |^2} \le 4\kappa _f^2,
		\end{align*}
		which satisfying the conditions 
		%	\begin{align*}
		$	\mathbb{E}[ {| {{X_i}} |} ] = 0, \mathbb{E}[ {X_i^2} 
		] \le 4\kappa _f^2.
		$
		%	\end{align*}
		Define the new variable
		\begin{align*}
		{Z_1} =& {f_s}(x) - f\left( x \right) = \frac{1}{{\left| {{\mathcal{S}_h}} 
				\right|}}\sum\nolimits_{j \in \left| {{\mathcal{S}_h}} \right|} {\left( 
			{{f_i}\left( x \right) - f\left( x \right)} \right)}, \\
		{Z_2} =& \sum\nolimits_{j \in \left| {{\mathcal{S}_h}} \right|} {\left( 
			{{f_i}\left( x \right) - f\left( x \right)} \right)}. 
		\end{align*}	
		Based on the Operator-Bernstein inequality \cite{gross2010note}, We 
		give  probability about the condition in Assumption	
		\ref{Newton:Assumption:approximation}, we have $
		| {{Z_2}} | = \left| {{\mathcal{S}_h}} \right|\left| {{Z_1}} 
		\right| \ge \left| {{\mathcal{S}_h}} \right|{\varepsilon _h}{\|s_k\|^2}$, then 
		%\begin{align*}
		$
		\Pr \left[ {\left| {{Z_2}} \right| > \left| {{\mathcal{S}_h}} 
			\right|{\varepsilon _h}{\Delta^2 _{{k}}}} \right] \le 2d\exp 
		\left( { - \frac{{{{\left( {\left| {{\mathcal{S}_h}} \right|{\varepsilon 
									_h}{\|s_k\|^2}} 
							\right)}^2}}}{{4\left| {{\mathcal{S}_h}} 
					\right|4\kappa _f^2}}} \right) \le \delta.
		$ 
		Thus, the cardinality of $\mathcal{S}_h$ should satisfy $\left| {{\mathcal{S}_h}} \right| 
		\ge 
		\frac{{16\kappa _f^2}}{{\varepsilon 
				_h^2\|s_k\|^4}}\log \left( {\frac{2d}{\delta }} 
		\right),$
		
		Furthermore,  based on Lemma \ref{Newton:Appendix:Tool:RandomSubset} and 
		Assumption 
		\ref{Newton:Assumption:Bound:Variance},  we can also have the upper 
		bounds with respect to the gradient 
		and the Hessian of $h(x)$,
		\[ 	{\left\| {\nabla h\left( x \right) - \nabla f\left( x \right)} 
			\right\|^2} \le 
		\frac{{\mathbb{I}\left( {\left| {{\mathcal{S}_h}} \right| < n} 
				\right)}}{{\left| 
				{{\mathcal{S}_h}} \right|}}{H_1},
		{\left\| {{\nabla ^2}h\left( x \right) - {\nabla ^2}f\left( x \right)} 
			\right\|^2} \le
		\frac{{\mathbb{I}\left( {\left| {{\mathcal{S}_h}} \right| < n} 
				\right)}}{{\left| 
				{{\mathcal{S}_h}} \right|}}{H_2}.\]	
	\end{proof}
	%%%%%%%%%%%%%%%%%%%%%%%%%%%%%%%%%%%%%%%%%%%%%%%%%%%%%%%%%%%%%%%%%
	%%%%%%%%%%%%%%%%%%%%%%%%%%%%%%%%%%%%%%%%%%%%%%%%%%%%%%%%%%%%%%%%%%%%%
	%%%%%%%%%%%%%%%%%%%%%%%%%%%%%%%%%%%%%%%%%%%%%%%%%%%%%%%%%%%%%%%%%%%
	%%%%%%%%%%%%%%%%%%%%%%%%%%%%%%%%%%%%%%%%%%%%%%%%%%%%%%%%%%%%%%%%%
	%%%%%%%%%%%%%%%%%%%%%%%%%%%%%%%%%%%%%%%%%%%%%%%%%%%%%%%%%%%%%%%%
	\section{Proof for Stochastic Trust Region}\label{Newton:Appendix:STR}
	%%%=============================================================================================
	%%%=============================================================================================
	%%%=============================================================================================
	\textbf{Proof of Lemma \ref{Newton:STR:Lemma:Upperbound_mx-msx}}
	\begin{proof}
		For the case $x_k\in \mathcal{S}_{\nabla f}$,
		through adding and subtracting the term $\nabla f( x_k ) $, we 
		have the lower bound of $\|g(x)\|$,
		\begin{align}
		\| {g( x_k )} \| &= \| {g( x_k ) -	 \nabla f( x_k ) + \nabla f( x_k )} 
		\|\nonumber\\
		&\ge \| {\nabla f( x_k )} \| -  \| {g( x_k 
			) - \nabla f( x_k )} \|\nonumber\\
		%\label{Newton:STR:Lemma:Upperbound_mx-msx:inequality1}
		&\ge {\varepsilon _{\nabla f}} - {\varepsilon_g},\nonumber
		\end{align}
		where the last inequality is based on the approximation of $\nabla f(x_k)$ in Assumption 
		\ref{Newton:Assumption:approximation}. Following the lower bound 
		on 	the decrease of the proximal quadratic 
		function ${m_k}( {{s}} )$ from  
		(4.20) in \cite{nocedal2006numerical}, we have
		\begin{align*}
		m_k(0) - {m_k}(s_k) 
		&\ge \frac{1}{2}\| {g( {{x_k}} )} 
		\|\min 	\left\{ {{\Delta _k},\frac{{\| {g( {{x_k}} )} 
					\|}}{{\| {B( {{x_k}} )} 
					\|}}} \right\}
		%\\& 
		\ge \frac{1}{2}({{ \varepsilon _{\nabla f}} - {\varepsilon_g})}\min 
		\left\{ {{\Delta _k},\frac{{{\varepsilon _{\nabla f}} - {\varepsilon_g
			}}}{{{\kappa _H}}}} \right\},
		\end{align*}
		where the last inequality is from above inequality and the bound 
		of 
		${\left\| {B({x_k})} \right\|}$ in Assumption 
		\ref{Newton:Assumption:bound}.
		
		For the case $x_k\in \mathcal{S}_H$,	through adding and subtracting the term $\nabla 
		^2f(x_k 
		) $, we 
		have
		\begin{align}
		\frac{{s_k^{T} B({x_k}){s_k}}}{{\left\| {{s_k}} \right\|}} =& 
		\frac{{s_k^T(B({x_k}) - \nabla^2f(x_k) + \nabla^2f(x_k)){s_k}}}{{\left\| {{s_k}} 
				\right\|}}\nonumber\\
		=& \frac{{s_k^T(B({x_k}) - \nabla^2f(x_k)){s_k}}}{{\left\| {{s_k}} \right\|}} + 
		\frac{{s_k^T\nabla^2f(x_k){s_k}}}{{\left\| {{s_k}} \right\|}}\nonumber\\
		\le& \left\| { B({x_k}) - \nabla^2f(x_k)} \right\| + 
		\frac{{s_k^T\nabla^2f(x_k){s_k}}}{{\left\| {{s_k}} \right\|}}\nonumber\\
		\le& {\varepsilon _B} - {v_0}\nonumber\\
		\label{Newton:STR:Lemma:Upperbound_fx-hsx_Lambda:inequality1}
		\le& {\varepsilon _B} - {\varepsilon _H} =  - \left( {{\varepsilon _H} 
			- 
			{\varepsilon _B}} \right),
		\end{align}
		where the second inequality is based on the Assumption in 
		\ref{Newton:Assumption:approximation} and Rayleigh quotient 
		\cite{conn2000trust}. Using the definition of 
		$m_k(s)$ in (\ref{Newton:Definition_m}), we have
		\begin{align*}
		m_k( 0) - {m_k}( {{s_k}} ) 
		%	&
		%	= f( 
		%	{{x_k}} ) - \left( {f( {{x_k}} ) + \langle {g({x_k}),{s_k}} \rangle  + 
		%		\frac{1}{2}s_k^TB( {{x_k}} ){s_k}} \right)\\
		&=  - \left\langle {g({x_k}),{s_k}} \right\rangle  - 
		\frac{1}{2}s_k^TB\left( {{x_k}} \right){s_k}\\
		&\mathop  \ge \limits^{{\scriptsize \textcircled{\tiny{1}}}}  - 
		\frac{1}{2}s_k^TB\left( {{x_k}} \right){s_k}\\
		&\mathop  \ge \limits^{{\scriptsize \textcircled{\tiny{2}}}} 
		\frac{1}{2}\left( {{\varepsilon _H} - {\varepsilon _B}} 
		\right){\left\| 
			{{s_k}} \right\|^2}
		%\\
		%&
		\mathop =\limits^{{\scriptsize \textcircled{\tiny{3}}}} \frac{1}{2} 
		\left( {{\varepsilon _H} - {\varepsilon _B}} 
		\right)\Delta _k^2,
		\end{align*}
		where inequality $\scriptsize \textcircled{\tiny{1}}$ is based on  $\langle 
		{g({x_k}),{s_k}} \rangle \le 0$, inequality $\scriptsize 
		\textcircled{\tiny{2}}$ follows from 
		(\ref{Newton:STR:Lemma:Upperbound_fx-hsx_Lambda:inequality1}), 
		inequality $\scriptsize 
		\textcircled{\tiny{3}}$ is based on $\left\| {{s_k}} \right\| = {\Delta 
			_k}$.
	\end{proof}
	%%%=============================================================================================
	%%%=============================================================================================
	%%%=============================================================================================
	\textbf{Proof of Lemma \ref{Newton:STR:Lemma:Upperbound_fxs-msx}}
	\begin{proof}
		Consider the Taylor expansion for $h(x_k+s_k)$ at $x_k$,
		\[h\left( {{x_k} + {s_k}} \right) = h\left( {{x_k}} \right) + \left\langle 
		{{s_k},\nabla h\left( {{x_k}} \right)} \right\rangle  + 
		\frac{1}{2}s_k^T\nabla^2h\left( {{\xi _k}} \right){s_k},\]
		where ${\xi  _k} \in \left[ {{x_k},{x_k} + {s_k}} \right]$.	
		Based on the definition of $m_k(s)$ in 
		(\ref{Newton:Definition_m}) and the Taylor expansion of the  
		function $ h(x_{k}+s_k)$ at $x_k$ above, we have 	
		\begin{align}
		&{m_k}\left( 0 \right) - {m_k}\left( {{s_k}} \right) - \left( {h\left( {{x_k}} \right) - 
			h\left( {{x_k} + {s_k}} \right)} \right)\nonumber\\
		=&  - \langle g({x_k}),{s_k}\rangle  - \frac{1}{2}s_k^TB({x_k}){s_k} + \left( {\left\langle 
			{\nabla h\left( {{x_k}} \right),{s_k}} \right\rangle  + \frac{1}{2}s_k^T{\nabla 
			^2}h\left( 
			{{\xi _k}} \right){s_k}} \right)\nonumber\\
		\label{Newton:STR:Lemma:Upperbound_fxs-msx-equality}
		=& \left\langle {\nabla h\left( {{x_k}} \right) - g({x_k}),{s_k}} \right\rangle  + 
		\frac{1}{2}s_k^T\left( {{\nabla ^2}h\left( {{\xi _k}} \right) - B({x_k})} 
		\right){s_k}\\
		=& \left\langle {\nabla h\left( {{x_k}} \right) - \nabla f\left( {{x_k}} \right) + \nabla 
			f\left( {{x_k}} \right) - g({x_k}),{s_k}} \right\rangle \nonumber\\
		&+ \frac{1}{2}s_k^T\left( {{\nabla ^2}h\left( {{\xi _k}} \right) - {\nabla ^2}f\left( {{\xi 
					_k}} \right) + {\nabla ^2}f\left( {{\xi _k}} \right) - {\nabla ^2}f\left( 
					{{x_k}} 
			\right) + 
			{\nabla ^2}f\left( {{x_k}} \right) - B({x_k})} \right){s_k}\nonumber\\
		\mathop  \le \limits^{{\scriptsize \textcircled{\tiny{1}}}}& 2\left( {\left\| {\nabla 
				h\left( {{x_k}} \right) - \nabla f\left( {{x_k}} \right)} 
			\right\| + \left\| {\nabla f\left( {{x_k}} \right) - g({x_k})} \right\|} \right)\left\| 
		{{s_k}} \right\|\nonumber\\
		&+ \frac{3}{2}\left( {\left\| {{\nabla ^2}h\left( {{\xi _k}} \right) - {\nabla ^2}f\left( 
				{{\xi _k}} \right)} \right\| + {\left\| {{\nabla ^2}f\left( {{\xi _k}} \right) - 
					{\nabla 
						^2}f\left( {{x_k}} \right)} \right\|} + \left\| {{\nabla ^2}f\left( 
				{{x_k}} \right) - B({x_k})} \right\|} \right){\left\| {{s_k}} \right\|^2}\nonumber\\
		\mathop  \le \limits^{{\scriptsize \textcircled{\tiny{2}}}}& 2\left( 
		{\frac{{\mathbb{I}\left( 
					{\left| {{\mathcal{S}_h}} \right| < n} \right)}}{{\left| {{\mathcal{S}_h}} 
					\right|}}H_1 + 
			{\varepsilon _g}} 
		\right)\left\| {{s_k}} \right\| + \frac{3}{2}\left( {\frac{{\mathbb{I}\left( {\left| 
						{{\mathcal{S}_h}} \right| 
						< n} \right)}}{{\left| {{\mathcal{S}_h}} \right|}}{H_2} + L_H\left\| 
						{{s_k}} 
			\right\| + 
			{\varepsilon 
				_B}} \right){\left\| {{s_k}} \right\|^2}\nonumber\\
		\mathop  \le \limits^{{\scriptsize \textcircled{\tiny{3}}}}& 2\left( 
		{\frac{{\mathbb{I}\left( 
					{\left| {{\mathcal{S}_h}} \right| < n} \right)}}{{\left| {{\mathcal{S}_h}} 
					\right|}}H_1 + 
			{\varepsilon _g}} 
		\right){\Delta _k} + \frac{3}{2}\left( {\frac{{\mathbb{I}\left( {\left| {{\mathcal{S}_h}} 
						\right| < 
						n} 
					\right)}}{{\left| {{\mathcal{S}_h}} \right|}}{H_2} + L_H{\Delta _k} + 
					{\varepsilon 
				_B}} 
		\right)\Delta 
		_k^2\nonumber,
		\end{align}
		where inequality $\scriptsize \textcircled{\tiny{1}}$ follows from the 
		Holder's inequality, 
		inequality $\scriptsize \textcircled{\tiny{2}}$  is based on the approximation of $\nabla 
		f(x)$ 
		and $\nabla^2f(x)$ in Assumption \ref{Newton:Assumption:approximation}, and the Lipschitz 
		continuous 
		of Hessian matrix of $f(x)$ in Assumption  
		\ref{Newton:Assumption:function}, inequality $\scriptsize 
		\textcircled{\tiny{3}}$ follows from the constraint condition as in the 
		objective (\ref{Newton:STR:Objective}).
		
		Furthermore, for the case of $\mathcal{S}_h=\mathcal{S}_g$, the first term of 
		(\ref{Newton:STR:Lemma:Upperbound_fxs-msx-equality}) is equal to zero, then we have
		\begin{align*}
		{m_k}\left( 0 \right) - {m_k}\left( {{s_k}} \right) - \left( {h\left( {{x_k}} \right) - 
			h\left( {{x_k} + {s_k}} \right)} \right) \le\frac{3}{2}\left( {\frac{{\mathbb{I}\left( 
					{\left| {{\mathcal{S}_h}} \right| < 
						n} 
					\right)}}{{\left| {{\mathcal{S}_h}} \right|}}{H_2} + {\Delta _k} + {\varepsilon 
				_B}} 
		\right)\Delta 
		_k^2.
		\end{align*}	
	\end{proof}
	%%%=============================================================================================
	%%%=============================================================================================
	%%%=============================================================================================
	\textbf{Proof of Lemma \ref{Newton:STR:Lemma:Bound:Radius}}
	\begin{proof}	
		By setting $\frac{{\mathbb{I}\left( {\left| {{\mathcal{S}_h}} \right| < n} 
		\right)}}{{\left| 
				{{\mathcal{S}_h}} 
				\right|}}{H_1} \le {\varepsilon _g}$, and $\frac{{\mathbb{I}\left( {\left| 
					{{\mathcal{S}_h}} 
					\right| < n} 
				\right)}}{{\left| 
				{{\mathcal{S}_h}} \right|}}{H_2} \le {\varepsilon _B}$,  we consider 
		two cases:
		
		For the case $\nabla f(x)>\varepsilon_{\nabla f(x)}$, 	we assume that 
		${\Delta _k} \le \frac{{{\varepsilon _{\nabla f}} - {\varepsilon 
					_g}}}{{{\kappa _H}}}$, which is used for Lemma 
		\ref{Newton:STR:Lemma:Upperbound_mx-msx}. Combine with (\ref{Newton:STR:1-p}) and
		the results in Lemma 
		\ref{Newton:STR:Lemma:Upperbound_mx-msx}, Lemma 
		\ref{Newton:STR:Lemma:Upperbound_fxs-msx}, and $\|s_k\|^2<\Delta_k^2$, 
		we have
		\begin{align}
		%		1 - \rho  =& \frac{{{m_k}\left( 0 \right) - {m_k}\left( {{s_k}} \right) - \left( 
		%{h\left( 
		%		{{x_k}} \right) - h\left( {{x_k} + {s_k}} \right)} \right) + 2{\varepsilon 
		%_h}{\Delta 
		%		_k}}}{{{m_k}\left( 0 \right) - {m_k}\left( {{s_k}} \right)}}\\
		%		\le& \frac{{2\left( {\frac{{\mathbb{I}\left( {\left| {{S_h}} \right| < n} 
		%\right)}}{{\left| 
		%		{{S_h}} 
		%		\right|}}{H_1} + {\varepsilon _g}} \right){\Delta _k} + \frac{3}{2}\left( 
		%		{\frac{{\mathbb{I}\left( 
		%		{\left| {{S_h}} \right| < n} \right)}}{{\left| {{S_h}} \right|}}{H_2} + L_H{\Delta 
		%_k} 
		%+ 
		%		{\varepsilon _B}} \right)\Delta _k^2 + 2{\varepsilon _h}{\Delta 
		%_k}}}{{\frac{1}{2}\left( 
		%		{{\varepsilon _{\nabla f}} - {\varepsilon _g}} \right){\Delta _k}}}\\
		%		=& \frac{{2\left( {2{\varepsilon _g} + {\varepsilon _h}} \right) + 
		%\frac{3}{2}\left( 
		%		{L_H{\Delta _k} + 2{\varepsilon _B}} \right){\Delta _k}}}{{\frac{1}{2}\left( 
		%{{\varepsilon 
		%		_{\nabla f}} - {\varepsilon _g}} \right)}}.
		1 - \rho 
		\le& \frac{{2\left( {\frac{{\mathbb{I}\left( {\left| {{S_h}} \right| < n} \right)}}{{\left| 
							{{S_h}} 
							\right|}}{H_1} + {\varepsilon _g}} \right){\Delta _k} + 
							\frac{3}{2}\left( 
				{\frac{{\mathbb{I}\left( 
							{\left| {{S_h}} \right| < n} \right)}}{{\left| {{S_h}} \right|}}{H_2} + 
					L_H{\Delta _k} + 
					{\varepsilon _B}} \right)\Delta _k^2 + 2{\varepsilon _h}{\Delta^2 
					_k}}}{{\frac{1}{2}\left( 
				{{\varepsilon _{\nabla f}} - {\varepsilon _g}} \right){\Delta _k}}}\nonumber\\
		\label{Newton:STR:Lemma:Bound:Radius-3}
		=&\frac{{4{\varepsilon _g} + \frac{3}{2}\left( {{L_H}{\Delta _k} + 2{\varepsilon _B} + 
					\frac{4}{3}{\varepsilon 
						_h}} \right){\Delta _k}}}{{\frac{1}{2}\left( {{\varepsilon _{\nabla f}} - 
					{\varepsilon _g}} 
				\right)}}.
		\end{align}
		In order to have the lower bound radius $\Delta_k$ such that $1 - \rho  \le1-\eta$, we 
		consider 
		the parameters' setting:
		\begin{itemize}
			\item For the first term $4{\varepsilon _g} = 
			\frac{1}{4}\left( {1 - \eta } \right)\left( {{\varepsilon _{\nabla f}} - {\varepsilon 
			_g}} 
			\right)$, then we define, 
			${\varepsilon _g} = \frac{1}{{16}}\left( {1 - \eta } \right)\left( {{\varepsilon 
			_{\nabla 
						f}} - {\varepsilon _g}} \right)$.
			\item For the second term $\frac{3}{2}\left( {L_H{\Delta _k} + 2{\varepsilon _B}+ 
				{\varepsilon _h}} 
			\right){\Delta 	_k} \le \frac{1}{4}\left( {1 - \eta } \right)\left( {{\varepsilon 
			_{\nabla 
						f}} - {\varepsilon _g}} \right)$, as $\varepsilon_B<1,\varepsilon 
			_h<1,\varepsilon_{\nabla f}-\varepsilon_g<1$,
			Thus, as long as 
			\begin{align}
			\label{Newton:STR:Lemma:Bound:Radius-1}
			{\Delta _k} \le \left( {{\varepsilon _{\nabla f}} - {\varepsilon _g}} \right)\min 
			\left\{ 
			{\frac{1}{{40}}\left( {1 - \eta } \right),\sqrt {\frac{1}{{12{L_H}}}\left( {1 - \eta } 
					\right)} } \right\}
			\end{align}
			we can obtain that $1 - \rho  \le1-\eta$. 	
		\end{itemize}
		At the $k$-iteration, when the radius $\Delta_k$ satisfy above condition, the update is 
		successful 
		iteration, as in Algorithm \ref{Newton:STR:Algorithm}, the radius $\Delta_k$ will increase 
		by 
		a factor $r_2$. 
		%However, 
		%if ${\Delta _k}>$(\ref{Newton:STR:Lemma:Bound:Radius-1}), it may 
		%or not be a successful iteration. If it is unsuccessful iteration, $\Delta_k$ will 
		%decrease by 
		%a 
		%factor $r_1$. Thus, in order to ensure the radius leading to successful iteration, 
		
		For the case ${\lambda _{\text{min} }}( {\nabla^2f(x_k)} ) \le  - 
		{\varepsilon _H}$, based on the results in Lemma 
		\ref{Newton:STR:Lemma:Upperbound_mx-msx} and Lemma 
		\ref{Newton:STR:Lemma:Upperbound_fxs-msx}, we have
		\begin{align*}
		1 - \rho  \le \frac{{\frac{3}{2}\left( {{L_H}{\Delta _k} + 2{\varepsilon _B} + 
					\frac{4}{3}{\varepsilon _h}} \right){\Delta _k}}}{{\frac{1}{2}\left( 
					{{\varepsilon 
						_H} 
					- {\varepsilon _B}} \right){\Delta _k}}}.
		\end{align*}
		In order to have the lower bound radius $\Delta_k$ such that $1 - \rho  \le1-\eta$, we 
		consider 
		the parameters' setting:
		\begin{align}\label{Newton:STR:Lemma:Bound:Radius-2}
		{\varepsilon _B} = {\varepsilon _h} = \frac{1}{{20}}\left( {1 - \eta } \right)\left( 
		{{\varepsilon _H} - {\varepsilon _B}} \right),{\Delta _k} \le \frac{1}{{6{L_H}}}\left( {1 - 
			\eta } \right)\left( {{\varepsilon _H} - {\varepsilon _B}} \right).
		\end{align}
		Based on above analysis and combine the assumption bound of $\Delta$ at beginning, there 
		exist 
		a 
		minimal radius 
		\[{\Delta _{{\rm{min}}}} = \min \left\{ {{\varepsilon _{\nabla f}} - {\varepsilon 
				_g},{\varepsilon 
				_H} - {\varepsilon _B}} \right\}{\kappa _1},{\kappa _1} = {r_1}\min \left\{ 
		{\frac{1}{{{\kappa 
						_H}}},\frac{1}{{40}}\left( {1 - \eta } \right),\sqrt 
						{\frac{1}{{12{L_H}}}\left( 
				{1 - \eta } 
				\right)} ,\frac{1}{{10{L_H}}}\left( {1 - \eta } \right)} \right\},\]
		where $0<r_1<1$. (multiply $r_1$ due to the fact that 
		(\ref{Newton:STR:Lemma:Bound:Radius-1}) 
		and 
		(\ref{Newton:STR:Lemma:Bound:Radius-2})
		plus a small constant may lead to a successful iteration such that   $\Delta_k$ will be 
		decreased the by 
		a factor $r_1$.)					
	\end{proof}
	%%%=============================================================================================
	%%%=============================================================================================
	%%%=============================================================================================
	\textbf{Proof of Theorem \ref{Newton:STR:Theorem:Iteration}}
	\begin{proof}
		Consider two index sets: $\mathcal{S}_{\nabla f}$ and $\mathcal{S}_{\lambda}$, we 
		separately 
		analyze the number of successful iteration based on results in Lemma 
		\ref{Newton:STR:Lemma:Upperbound_mx-msx} and
		\ref{Newton:STR:Lemma:Bound:Radius}. And then add both of them to form 
		the most number of successful iterations. Let  $f_\text{low}$ is the minimal value of 
		objective, we 
		have two kinds of successful iterations:	
		\begin{itemize}
			\item Consider the case of $\|	{\nabla f( x )} \| \ge {\varepsilon 
				_{\nabla f}}$, if $k$th iteration is successful, then we have
			\begin{align*}
			f\left( {{x_k}} \right) - f\left( {{x_k} + {s_k}} \right) &\ge {\eta 
			}\left( {{m_k}\left(0 \right) - {m_k}\left( { 
					{s_k}} 
				\right)} \right)\\
			&\ge \frac{1}{2}\eta ({\varepsilon _{\nabla f}} - {\varepsilon _g}){
			}\min \left\{ {{\Delta _k},\frac{{{\varepsilon _{\nabla f}} - 
						{\varepsilon _g}}}{{{\kappa _H}}}} \right\}
			\\&
			\ge\frac{1}{2}{\eta}({\varepsilon _{\nabla f}} - {\varepsilon _g}){\Delta 
			_{\text{min1}}},
			\end{align*}
			where the last two inequalities are based on Lemma 
			\ref{Newton:STR:Lemma:Upperbound_mx-msx} and Lemma 
			\ref{Newton:STR:Lemma:Bound:Radius}. Let $T_{1}$ denotes 
			the number of successful iterations for $k\in \mathcal{S}_{\nabla f}$. Applying above 
			inequality, we can obtain 
			\begin{align*}
			f\left( {{x_0}} \right) - {f_\text{low}}
			%		 &\ge \sum\limits_{k = 0}^{{T_1}} 
			%		{\left( {f\left( {{x_k}} \right) - f\left( {{x_k} + {s_k}} \right)} 
			%			\right)} 
			%\\&
			\ge \frac{1}{2}{T_1}({\varepsilon _{\nabla f}} - {\varepsilon _g}){\eta}{\Delta 
				_\text{min1}}.
			\end{align*}
			where $f_\text{low}$ is blower bound of the objective.
			\item Consider the case of ${\lambda _{\text{min} }}( \nabla^2f(x_k) ) 
			\le  - 	{\varepsilon _H}$, if $k$th iteration is successful, based on Lemma 
			\ref{Newton:STR:Lemma:Upperbound_mx-msx} and Lemma 
			\ref{Newton:STR:Lemma:Bound:Radius}, we have
			\begin{align*}
			f\left( {{x_k}} \right) - f\left( {{x_k} + {s_k}} \right) &\ge {\eta 
			}\left( {{m_k}\left(0 \right) - {m_k}\left( { 
					{s_k}} 
				\right)} \right)
			%\\&
			\ge \frac{1}{2}\eta({\varepsilon _H} - {\varepsilon _B})\Delta _\text{min2}^2.
			\end{align*}
			Let $T_2$  denote 
			the number of successful iteration for $k\in \mathcal{S}_{H}$, we obtain 
			\begin{align*}
			f\left( {{x_0}} \right) - {f_{low}}
			%		 &\ge \sum\limits_{k = 
			%			0}^{{T_2}} {\left( {f\left( {{x_k}} \right) - f\left( {{x_k} + 
			%					{s_k}} \right)} \right)} 
			%\\& 
			\ge \frac{1}{2}{T_2}{\eta}({\varepsilon _H} - {\varepsilon 
				_B})\Delta _{\text{min2}}^2.
			\end{align*}
		\end{itemize}
		Let $T_\text{suc}$ denotes the number of  successful iterations, combing 
		above 
		iteration and $\Delta_{\text{min}}$ in Lemma \ref{Newton:STR:Lemma:Bound:Radius},  we have
		\begin{align*}
		T_\text{suc} \le {T_1} + {T_2} \le& \frac{{2\left( {f\left( {{x_0}} \right) - 
					{f_{{\rm{low}}}}} \right)}}{{({\varepsilon _{\nabla f}} - {\varepsilon _g})\eta 
				{\Delta 
					_\text{min1}}}} + \frac{{2\left( {f\left( {{x_0}} \right) - {f_{{\rm{low}}}}} 
				\right)}}{{\eta 
				({\varepsilon _H} - {\varepsilon _B}){\Delta^2_\text{min2}}}}\\
		\le& \frac{{2\left( {f\left( {{x_0}} \right) - {f_{{\rm{low}}}}} 
				\right)}}{{({\varepsilon _{\nabla f}} - {\varepsilon _g})\eta {\kappa _1}\left( 
				{{\varepsilon 
						_{\nabla f}} - {\varepsilon _g}} \right)}} + \frac{{2\left( {f\left( 
						{{x_0}} 
					\right) - 
					{f_{{\rm{low}}}}} \right)}}{{({\varepsilon _H} - {\varepsilon _B})\eta {{\left( 
						{{\kappa _2}\left( 
							{{\varepsilon _H} - {\varepsilon _B}} \right)} \right)}^2}}}\\
		\le& {\kappa _3}{\rm{max}}\left\{ {{{\left( {{\varepsilon _{\nabla f}} - {\varepsilon _g}} 
					\right)}^{ - 2}},{{({\varepsilon _H} - {\varepsilon _B})}^{ - 3}}} \right\},
		\end{align*}
		where ${\kappa _3} = 2\left( {f\left( {{x_0}} \right) - {f_{low}}} \right){\rm{max}}\left\{ 
		{1/\left( {\eta {\kappa _1}} \right),1/\left( {\eta \kappa _2^2} \right)} \right\}$.

		Let $T_\text{unsuc}$ denotes the number of unsuccessful iteration, we have $	{{ 
				r}_1}{\Delta 
			_k} \le 
		{\Delta _{k + 1}}$;	Let $T_\text{suc}$ denotes the number of successful iteration, we have 
		${r_2}{\Delta 
			_k} \le {\Delta _{k + 1}}$.	Thus, we inductively deduce,
		\[{\Delta _{\text{min} }}r_2^{T_\text{suc}}r_1^{T_\text{unsuc}} \le {\Delta _{\text{max} 
		}}.\]
		where $\Delta _\text{max}$ is defined in  (\ref{Newton:STR:Lemma:Bound:Radius-equality}) 
		and 
		${\Delta _{\text{min} }}$ is defined in Algorithm \ref{Newton:STR:Algorithm}.
		Thus, the number of unsuccessful index set is at most
		\[{T_\text{unsuc}} \ge \frac{1}{{ - \log {r_1}}}\left( {\log \left( {\frac{{{\Delta 
							_{{\text{max}}}}}}{{{\Delta 
							_{{\text{min}}}}}}} \right) - {T_\text{suc}}\log {r_2}} \right).\]
		Combine with the successful iteration, we can obtain the total 
		iteration complexity,
		\begin{align*}
		{T_\text{suc}} + {T_\text{unsuc}} =& {T_\text{suc}}\left( {1 + \log \left( {\frac{{{\Delta 
							_{{\rm{max}}}}}}{{{\Delta 
							_{{\rm{min}}}}}}} \right)} \right) - \log \left( {\frac{{{\Delta 
						_{{\rm{max}}}}}}{{{\Delta 
						_{{\rm{min}}}}}}} \right)\frac{1}{{\log {r_1}}}\\
		=& O\left( {\max \left\{ {{{({\varepsilon _{\nabla f}} - {\varepsilon _g})}^{ - 
						2}},{{({\varepsilon _H} - {\varepsilon _B})}^{ - 3}}} \right\}} \right).
		\end{align*}
	\end{proof}
	%%%%%%%%%%%%%%%%%%%%%%%%%%%%%%%%%%%%%%%%%%%%%%%%%%%%%%%%%%%%%%%%%
	%%%%%%%%%%%%%%%%%%%%%%%%%%%%%%%%%%%%%%%%%%%%%%%%%%%%%%%%%%%%%%%%%%%%%
	%%%%%%%%%%%%%%%%%%%%%%%%%%%%%%%%%%%%%%%%%%%%%%%%%%%%%%%%%%%%%%%%%%%
	%%%%%%%%%%%%%%%%%%%%%%%%%%%%%%%%%%%%%%%%%%%%%%%%%%%%%%%%%%%%%%%%%
	%%%%%%%%%%%%%%%%%%%%%%%%%%%%%%%%%%%%%%%%%%%%%%%%%%%%%%%%%%%%%%%%

	%%%%%%%%%%%%%%%%%%%%%%%%%%%%%%%%%%%%%%%%%%%%%%%%%%%%%%%%%%%%%%%%%
	%%%%%%%%%%%%%%%%%%%%%%%%%%%%%%%%%%%%%%%%%%%%%%%%%%%%%%%%%%%%%%%%%%%%%
	%%%%%%%%%%%%%%%%%%%%%%%%%%%%%%%%%%%%%%%%%%%%%%%%%%%%%%%%%%%%%%%%%%%
	%%%%%%%%%%%%%%%%%%%%%%%%%%%%%%%%%%%%%%%%%%%%%%%%%%%%%%%%%%%%%%%%%
	%%%%%%%%%%%%%%%%%%%%%%%%%%%%%%%%%%%%%%%%%%%%%%%%%%%%%%%%%%%%%%%%
	\section{Proof for Stochastic ARC}\label{SFSNewton:Appendix:SARC}
	%%%===================================================
	%%%=====================================================
	%%%====================================================
	\textbf{Proof of Lemma \ref{Newton:SARC:Lemma:UpperboundOfP0Ps} }
	\begin{proof}
		Based on the Cauchy-Schwarz inequality and the definition of $m_k(x)$,	we 	have
		\begin{align*}
		{p_k}\left( {s_k^C} \right) \le& p_k\left( 0 \right) + 
		\left\langle 
		{g\left( {{x_k}} \right), - \alpha g\left( {{x_k}} \right)} 
		\right\rangle  
		+ \frac{1}{2}\left\langle {B\left( {{x_k}} \right)\left( { - \alpha 
				g\left( 
				{{x_k}} \right)} \right), - \alpha g\left( {{x_k}} \right)} 
		\right\rangle  
		+ \frac{{{\sigma _k}}}{3}{\left\| { - \alpha g\left( {{x_k}} \right)} 
			\right\|^3}\\
		\le& p_k\left( 0 \right) - {\left\| {g\left( {{x_k}} \right)} 
			\right\|^2}\alpha  + \frac{1}{2}\left\| {B\left( {{x_k}} \right)} 
		\right\|{\left\| {g\left( {{x_k}} \right)} \right\|^2}{\alpha ^2} + 
		\frac{{{\sigma _k}}}{3}{\left\| {g\left( {{x_k}} \right)} 
			\right\|^3}{\alpha ^3}\\
		=&p_k(0)+q(\alpha),
		\end{align*}
		where 
		\begin{align*}
		q\left( \alpha  \right) =&- {\left\| {g\left( {{x_k}} \right)} 
			\right\|^2}\alpha  + \frac{1}{2}\left\| {B\left( {{x_k}} \right)} 
		\right\|{\left\| {g\left( {{x_k}} \right)} \right\|^2}{\alpha ^2} + 
		\frac{{{\sigma _k}}}{3}{\left\| {g\left( {{x_k}} \right)} 
			\right\|^3}{\alpha ^3}\\
		=&  - {\left\| {{g_k}} \right\|^2}\alpha  + 
		\frac{1}{2}\left\| {{B_k}} \right\|{\left\| {{g_k}} \right\|^2}{\alpha 
			^2} + \frac{{{\sigma _k}}}{3}{\left\| {{g_k}} \right\|^3}{\alpha 
			^3}\\
		=& {\left\| {{g_k}} \right\|^2}\left( { - \alpha  + \frac{1}{2}\left\| 
			{{B_k}} \right\|{\alpha ^2} + \frac{{{\sigma _k}}}{3}\left\| 
			{{g_k}} 
			\right\|{\alpha ^3}} \right).
		\end{align*}
		For simplicity, we use ${g_k} = 
		g\left( {{x_k}} \right)$ and ${B_k} = B\left( {{x_k}} \right)$ instead.
		
		In order to show  ${p_k}\left( {s_k^C} \right) \le f\left( {{x_k}} 
		\right)$, we should to check that the \textbf{minimal value} of 
		$q(\alpha)$ 
		is 	negative or not. That is, if the minimal value of $q(\alpha)$ is 
		negative, there is 	exist a $\alpha_0$ such that $q(\alpha_0)\le0$, 
		and  
		${p_k}\left( {s_k^C}\right) \le f(x_k)+q(\alpha)$ for $\forall 
		\alpha>0$, 
		then, we can obtain ${p_k}\left( {s_k^C} \right) \le f\left( {{x_k}} 
		\right)$. 	Now, consider the gradient of $p(\alpha)$,
		\begin{align*}
		\nabla q\left( \alpha  \right) = {\left\| {{g_k}} \right\|^2}\left( { - 
			1 + \left\| {{B_k}} \right\|\alpha  + {\sigma _k}\left\| {{g_k}} 
			\right\|{\alpha ^2}} \right).
		\end{align*}
		Let $\nabla q\left( \alpha  \right) = 0$, we have two solutions,
		\begin{center}
			$	{\alpha _1} = \frac{{ - \left\| {{B_k}} \right\| - \sqrt {{{\left\| 
								{{B_k}} \right\|}^2} + 4{\sigma _k}\left\| {{g_k}} 
						\right\|} 
			}}{{2{\sigma _k}\left\| {{g_k}} \right\|}} < 0,{\alpha _2} = \frac{{ - 
					\left\| {{B_k}} \right\| + \sqrt {{{\left\| {{B_k}} 
								\right\|}^2} + 
						4{\sigma _k}\left\| {{g_k}} \right\|} }}{{2{\sigma 
						_k}\left\| {{g_k}} 
					\right\|}}.$
		\end{center}
		Because we require $\alpha >0$, we do not need to consider 
		$\alpha_1<0$. 
		Thus, we 
		obtain the geometrical character,
		\begin{center}
			$	\nabla q\left( \alpha  \right) = \left\{ {\begin{array}{*{20}{c}}
				+ &{\alpha  > {\alpha _2}}\\
				0&{\alpha  = {\alpha _2}}\\
				- &{{\alpha _2} > \alpha  \ge 0}
				\end{array}} \right. \Rightarrow q\left( \alpha  \right)=\left\{ 
			{\begin{array}{*{20}{c}}
				\nearrow &{\alpha  \ge {\alpha _2}}\\
				\searrow &{{\alpha _2} > \alpha  \ge 0}
				\end{array}} \right.,$
		\end{center}
		where we use $``+"$ and $``-"$ to denote the positive and negative of 
		$\nabla 
		q\left( \alpha  \right)$, respectively, use $``\nearrow "$ and $ 
		``\searrow" $ to 
		denote the increasing and decreasing functions, respective. From above 
		description, we can obtain that $\alpha_2$ is the minimal solution. 
		Putting 
		$\alpha_2$ into $q(\alpha)$, we have
		\begin{align*}
		q\left( {{\alpha _2}} \right) =& {\left\| {{g_k}} \right\|^2}\left( { - 
			{\alpha 
				_2} + \frac{1}{2}\left\| {{B_k}} \right\|\alpha _2^2 + 
			\frac{1}{3}{\sigma 
				_k}\left\| {{g_k}} \right\|\alpha _2^3} \right)\\
		=& \left\| {{g_k}} \right\|\left( {{\alpha _2}\left\| {{g_k}} 
			\right\|\left( { 
				- 
				1 + \left\| {{B_k}} \right\|{\alpha _2} + {\sigma _k}\left\| 
				{{g_k}} 
				\right\|\alpha _2^2} \right) - \frac{1}{2}\left\| {{B_k}} 
			\right\|\left\| 
			{{g_k}} \right\|\alpha _2^2 - \frac{2}{3}{\sigma _k}{{\left\| 
					{{g_k}} 
					\right\|}^2}\alpha _2^3} \right)\\
		=&  - {\left\| {{g_k}} \right\|^2}\alpha _2^2\left( {\frac{1}{2}\left\| 
			{{B_k}} 
			\right\| + \frac{2}{3}{\sigma _k}\left\| {{g_k}} \right\|{\alpha 
				_2}} \right)\\
		=&  - {\left\| {{g_k}} \right\|^2}\alpha _2^2\left( {\frac{1}{6}\left\| 
			{{B_k}} 
			\right\| + \frac{1}{3}\sqrt {{{\left\| {{B_k}} \right\|}^2} + 
				4{\sigma 
					_k}\left\| {{g_k}} \right\|} } \right)\\
		%	=&  - {\left\| {{g_k}} \right\|^2}\alpha _2^2\left( {\frac{1}{6}\left\| 
		%		{{B_k}} 
		%		\right\| + \frac{1}{6}\sqrt {{{\left\| {{B_k}} \right\|}^2} + 
		%			4{\sigma 
		%				_k}\left\| {{g_k}} \right\|}  + \frac{1}{6}\sqrt {{{\left\| 
		%					{{B_k}} 
		%					\right\|}^2} + 4{\sigma _k}\left\| {{g_k}} \right\|} } 
		%	\right)\\
		\le&  - \frac{1}{6}{\left\| {{g_k}} \right\|^2}\alpha _2^2\left( 
		{\left\| 
			{{B_k}} \right\| + \sqrt {{{\left\| {{B_k}} \right\|}^2} + 4{\sigma 
					_k}\left\| 
				{{g_k}} \right\|} } \right),
		\end{align*}
		where 	the third equality is from
		\begin{align*}
		{\left\| {{g_k}} \right\|\left( { - 1 + \left\| {{B_k}} \right\|{\alpha 
					_2} 
				+ {\sigma _k}\left\| {{g_k}} \right\|\alpha _2^2} \right)}=0.
		\end{align*}
		Because $\alpha_2$ can also be expressed as 
		\begin{center}
			${\alpha _2} = \frac{{ - \left\| {{B_k}} \right\| + \sqrt {{{\left\| 
								{{B_k}} 
								\right\|}^2} + 4{\sigma _k}\left\| {{g_k}} 
						\right\|} }}{{2{\sigma 
						_k}\left\| {{g_k}} \right\|}} = \frac{2}{{\left\| {{B_k}} 
					\right\| + \sqrt 
					{{{\left\| {{B_k}} \right\|}^2} + 4{\sigma _k}\left\| {{g_k}} 
						\right\|} }},$
		\end{center}
		we can obtain
		\begin{align*}
		q\left( {{\alpha _2}} \right) \le  - \frac{1}{3}{\left\| {{g_k}} 
			\right\|^2}{\alpha _2}
		=&  - \frac{1}{6}{\left\| {{g_k}} \right\|^2}\frac{2}{{\left\| {{B_k}} 
				\right\| + \sqrt {{{\left\| {{B_k}} \right\|}^2} + 4{\sigma 
						_k}\left\| 
					{{g_k}} \right\|} }}\\
		\le &- 
		\frac{1}{10}\left\| {{g_k}} \right\|\min \left\{ 
		{\frac{\left\| {{g_k}} \right\|}{{\left\| {{B_k}} 
					\right\|}},\sqrt {\frac{{\left\| {{g_k}} 
						\right\|}}{{{\sigma _k}}}} } 
		\right\},
		\end{align*}
		where inequality  is based on 
		difference 
		value between ${\left\| {{B_k}} \right\|^2}$ and ${\sigma _k}\left\| 
		{{g_k}} 
		\right\|$, and $\frac{1}{{3\left( {1 + \sqrt 5 } \right)}} 
		\approx 
		0.103 \ge \frac{1}{10}$.
		
		Furthermore, we can also see that 	Based on Lemma \ref{Newton:SARC:Lemma:UpperboundOfP0Ps} 
		and 
		the definition of $p_k(x)$ in 
		(\ref{Newton:SARC:definition_P}), we have
		\begin{align*}
		{p_k}\left(s_k \right) - p\left( 0 \right) = \left\langle 
		{s_k,g\left( {{x_k}} \right)} \right\rangle  + 
		\frac{1}{2}{s_k^T}B\left( 
		{{x_k}} \right)s_k + \frac{1}{3}{\sigma _k}{\left\| s_k \right\|^3} \le 
		0.
		\end{align*}
		Using Cauchy-Schwarz inequality, we have
		\begin{align*}
		\frac{1}{3}{\sigma _k}{\left\| s_k \right\|^3} \le  - \left\langle 
		{s_k,g\left( 
			{{x_k}} \right)} \right\rangle  - \frac{1}{2}{s_k^T}B\left( {{x_k}} 
		\right)s_k 
		\le 
		\left\| s_k \right\|\left\| {g\left( {{x_k}} \right)} \right\| + 
		\frac{1}{2}\left\| {B\left( {{x_k}} \right)} \right\|{\left\| s_k 
			\right\|^2}.
		\end{align*}
		Arranging the position of $\left\| s_k \right\|$, we have
		\begin{align}\label{Newton:SARC:Lemma:UpperboundOfP0Ps_Inequality}
		\left\| s_k \right\|( {\frac{1}{3}{\sigma _k}{{\left\| s_k 
					\right\|}^2} - 
			\frac{1}{2}\left\| {B\left( {{x_k}} \right)} \right\|\left\| s_k 
			\right\| - 
			\left\| {g\left( {{x_k}} \right)} \right\|} ) \le 0.
		\end{align}
		Because $\left\| s_k \right\|$ is positive, in order to satisfy above 
		equality, 
		$\left\| s_k \right\|$ is upper bounded by,
		%	\begin{center}
		%		$	\| s_k \| \le {\rm{ }}\frac{{3\| {B( {{x_k}} )} 
		%				\| + \sqrt {9{{\| {B( {{x_k}} )} \|}^2} + 48{\sigma 
		%						_k}\| {g( {{x_k}} )} \|} }}{{4{\sigma _k}}} \le 
		%		\frac{{11}}{4}\max \left\{ 
		%		{\frac{{\left\| {B({x_k})} \right\|}}{{{\sigma _k}}},\sqrt 
		%			{\frac{{\left\| 
		%						{g({x_k})} \right\|}}{{{\sigma _k}}}} } \right\},$
		%	\end{center}
		\begin{center}
			$\| s_k \| \le {\rm{ }}\frac{{3\| {B( {{x_k}} )} 
					\| + \sqrt {9{{\| {B( {{x_k}} )} \|}^2} + 48{\sigma 
							_k}\| {g( {{x_k}} )} \|} }}{{4{\sigma _k}}} \le 
			\frac{{11}}{4}\max \left\{ 
			{\frac{{\left\| {B({x_k})} \right\|}}{{{\sigma _k}}},\sqrt 
				{\frac{{\left\| 
							{g({x_k})} \right\|}}{{{\sigma _k}}}} } \right\},$
		\end{center}
		where the first inequality follows from the solution of a quadratic 
		function in
		(\ref{Newton:SARC:Lemma:UpperboundOfP0Ps_Inequality}), and the second 
		inequality is 
		based on the values between ${\left\| {B\left( {{x_k}} \right)} 
			\right\|}$ and 
		${{\sigma _k}\left\| {g\left( {{x_k}} \right)} \right\|}$.
	\end{proof}

	\textbf{Proof of Lemma \ref{Newton:SARC:Lemma:UpperboundOfP0Ps-2}}
	\begin{proof}
		Based on the definition of $p_k(s)$ in (\ref{Newton:SARC:definition_P}), we have
		\begin{align*}
		p( 0) - {p_k}( s_k ) &=  - 
		\left\langle 
		{g\left( {{x_k}} \right),{s_k}} \right\rangle  - 
		\frac{1}{2}s_k^TB\left( {{x_k}} \right){s_k} - \frac{{{\sigma 
					_k}}}{3}{\left\| {{s_k}} \right\|^3}\\
		&= \underbrace { - \left\langle {g\left( {{x_k}} \right),{s_k}} 
			\right\rangle  - s_k^TB\left( {{x_k}} \right){s_k} - {{\left\| 
					{{s_k}} 
					\right\|}^3}}_{ = 0} + \frac{1}{2}s_k^TB\left( {{x_k}} 
		\right){s_k} + 	\frac{{2{\sigma _k}}}{3}{\left\| {{s_k}} 
			\right\|^3}\\
		&\mathop  = \limits^{{\scriptsize 	\textcircled{\tiny{1}}}} 
		\frac{1}{2}s_k^TB\left( {{x_k}} \right){s_k} + 
		\frac{{2{\sigma _k}}}{3}{\left\| {{s_k}} \right\|^3}\\
		&\mathop  \ge \limits^{{\scriptsize \textcircled{\tiny{2}}}} 
		\frac{{{\sigma 
					_k}}}{6}{\left\| {{s_k}} 
			\right\|^3},
		\end{align*}
		where ${\scriptsize \textcircled{\tiny{1}}}$ and ${\scriptsize 
			\textcircled{\tiny{2}}}$ follows from 
		(\ref{Newton:SARC:Assumption_s1}) and 
		(\ref{Newton:SARC:Assumption_s2}), respectively.
		
		Consider the lower bound of $\|s_k\|$: Firstly, based on the definition of $m_k(s)$, for 
		simplicity, we use ${g_k} = 
		g\left( {{x_k}} \right)$ and ${B_k} = B\left( {{x_k}} \right)$ instead, 
		we 	have
		\begin{align*}
		\left\| {g\left( {{x_{k + 1}}} \right)} \right\| =& \left\| {g\left( 
			{{x_{k + 1}}} \right) 
			- \nabla {p_k}\left( s \right) + \nabla {p_k}\left( s \right)} 
		\right\|\\
		\mathop\le\limits^{{\scriptsize \textcircled{\tiny{1}}}}& \left\| 
		{g\left( {{x_k}} \right) 
			+ \int_0^1 {B\left( {{x_k} + \tau {s_k}} 
				\right){s_k}d\tau }  - \left( {g\left( {{x_k}} \right) + B\left( 
				{{x_k}} \right){s_k} + 
				{\sigma _k}\left\| {{s_k}} \right\|{s_k}} \right)} \right\| + \left\| 
		{\nabla {p_k}\left( 
			{{s_k}} \right)} \right\|\\
		\mathop\le\limits^{{\scriptsize \textcircled{\tiny{2}}}}& \left\| 
		{\int_0^1 {\left( 
				{B\left( {{x_k} + \tau {s_k}} \right) - H\left( {{x_k} + 
						\tau {s_k}} \right)} \right){s_k}d\tau } } \right\| + \left\| 
		{H\left( {{x_k} + \tau {s_k}} 
			\right)s - H\left( {{x_k}} \right){s_k}} \right\|\\
		&+ \left\| {H\left( {{x_k}} \right){s_k} - {B_k}{s_k}} \right\| + 
		{\sigma _k}{\left\| 
			{{s_k}} \right\|^2} + \left\| {\nabla {p_k}\left( {{s_k}} \right)} 
		\right\|\\
		\mathop\le\limits^{{\scriptsize \textcircled{\tiny{2}}}}& 
		{\varepsilon _B}\left\| {{s_k}} 
		\right\| + {L_H}{\left\| {{s_k}} \right\|^2} + 
		{\varepsilon _B}\left\| {{s_k}} \right\| + {\sigma _k}{\left\| 
			{{s_k}} \right\|^2} + 
		{\theta _k}\left\| {g\left( {{x_k}} \right)} \right\|\\
		=& 2{\varepsilon _B}\left\| {{s_k}} \right\| + \left( {{L_H} + 
			{\sigma _k}} \right){\left\| 
			{{s_k}} \right\|^2} + {\theta _k}\left\| {g\left( {{x_k}} \right)} 
		\right\|,
		\end{align*}	
		where inequality ${\scriptsize 
			\textcircled{\tiny{1}}}$ is based on the triangle inequality and  the Taylor expansion 
		of 
		$g(x)$,
		equality 
		${\scriptsize 	\textcircled{\tiny{2}}}$ is obtained by adding and 
		subtracting the term of 
		${H\left( {{x_k}} \right){s_k}}$ and ${H\left( {{x_k+\tau s_k}} 
			\right){s_k}}$, and triangle inequality; equality 
		${\scriptsize 	\textcircled{\tiny{3}}}$ follows from Assumption 
		\ref{Newton:Assumption:function}, \ref{Newton:Assumption:approximation} and the condition 
		in 
		(\ref{Newton:SARC:Assumption_s3}). Secondly, consider $g(x_k)$, we have
		\begin{align*}
		&\left\| {g\left( {{x_k}} \right)} \right\|\\
		=& \left\| {g\left( {{x_k}} \right) - \nabla f\left( {{x_k}} \right) + 
			\nabla f\left( {{x_k}} 
			\right) - \nabla f\left( {{x_k} + {s_k}} \right) + \nabla f\left( {{x_k} 
				+ {s_k}} \right) - 
			g\left( {{x_{k + 1}}} \right) + g\left( {{x_{k + 1}}} \right)} \right\|\\
		\le& \left\| {g\left( {{x_k}} \right) - \nabla f\left( {{x_k}} \right)} 
		\right\| + \left\| 
		{\nabla f\left( {{x_k}} \right) - \nabla f\left( {{x_k} + {s_k}} \right)} 
		\right\| + \left\| 
		{\nabla f\left( {{x_k} + {s_k}} \right) - g\left( {{x_{k + 1}}} \right)} 
		\right\| + \left\| 
		{g\left( {{x_{k + 1}}} \right)} \right\|\\
		\le& 2{\varepsilon _g} + {L_{\nabla f}}\left\| {{s_k}} \right\| + \left\| 
		{g\left( {{x_{k + 
						1}}} 
			\right)} \right\|,
		\end{align*}
		where the first and second inequality are based on the triangle 
		inequality, 
		Lipschitz continuity of gradient ${\nabla f\left( {{x_k}} \right)}$ in Assumption
		\ref{Newton:Assumption:function} and 
		\ref{Newton:Assumption:approximation}.
		
		Finally, replace the term $\left\| {\nabla g\left( {{x_k}} \right)} \right\|$, 
		we have
		\begin{align*}
		\left( {1 - {\theta _k}} \right)\left\| {g\left( {{x_{k + 1}}} \right)} 
		\right\| \le 
		2{\varepsilon _B}\left\| {{s_k}} \right\| + \left( {{L_H} + {\sigma _k}} 
		\right){\left\| 
			{{s_k}} \right\|^2} + 2{\theta _k}{\varepsilon _g} + {\theta 
			_k}{L_{\nabla f}}\left\| {{s_k}} 
		\right\|.
		\end{align*}
		Consider the definition of ${\theta _k} 
		\le 
		{\zeta _\theta }\min \left\{ {1,\left\| {{s_k}} \right\|} 
		\right\},{\zeta 
			_\theta } < 1$, we analysis the bound from different rang of 
		$\|s_k\|$
		\begin{itemize}
			\item For the case of $\left\| {{s_k}} \right\| \ge 1$, we have
			\begin{align*}
			\left( {1 - {\theta _k}} \right)\left\| {g\left( {{x_{k + 1}}} 
				\right)} \right\| \le 
			2{\varepsilon _B}{\left\| {{s_k}} \right\|^2} + \left( {{L_H} + 
				{\sigma _k}} 
			\right){\left\| {{s_k}} \right\|^2} + 2{\kappa _\theta }{\left\| 
				{{s_k}} 
				\right\|^2}{\varepsilon _g} + {\kappa _\theta }{L_{\nabla f}}{\left\| 
				{{s_k}} \right\|^2}.
			\end{align*}
			\item For the case of $\left\| {{s_k}} \right\| \le 1$, based 
			on the 
			assumption on 
			${\varepsilon_g}$ and ${\varepsilon _g}$, that is 
			\begin{align*}
			{\varepsilon _B} \le& {\zeta _1}\left( {{\varepsilon _{\nabla f}} - 
				{\varepsilon _g}} 
			\right) \le {\zeta _1}\left\| {g\left( {{x_{k + 1}}} \right)} 
			\right\|,\\
			{\varepsilon _g} \le& {\zeta _2}\left( {{\varepsilon _{\nabla f}} - 
				{\varepsilon _g}} 
			\right) \le {\zeta _2}\left\| {g\left( {{x_{k + 1}}} \right)} \right\|,
			\end{align*}
			where ${\zeta _1},{\zeta _2} < 1$. We have
			\begin{align*}
			\left( {1 - {\theta _k}} \right)\left\| {g\left( {{x_{k + 1}}} 
				\right)} \right\| \le& 
			2{\varepsilon _B} + \left( {{L_H} + {\sigma _k}} \right){\left\| 
				{{s_k}} \right\|^2} + 
			2{\kappa _\theta }{\varepsilon _g} + {\kappa _\theta }{L_{\nabla 
					f}}{\left\| {{s_k}} 
				\right\|^2}\\
			\le& \left( {2{\zeta _1} + 2{\kappa _\theta }{\zeta _2}} 
			\right)\left\| {g\left( {{x_{k + 
							1}}} \right)} \right\| + \left( {{L_H} + {\sigma _k} + {\kappa 
					_\theta }{L_{\nabla f}}} 
			\right){\left\| {{s_k}} \right\|^2}.
			\end{align*}
		\end{itemize}	
		Thus, in all, we can obtain
		\[\left\| {g\left( {{x_{k + 1}}} \right)} \right\| \le {\kappa 
			_s}{\left\| {{s_k}} 
			\right\|^2},{\kappa _s} = \min \left\{ {\frac{{2{\varepsilon _B} + \left( 
					{{L_H} + {\sigma _k}} 
					\right) + 2{\kappa _\theta }{\varepsilon _g} + {\kappa _\theta 
					}{L_{\nabla f}}}}{{\left( {1 - 
						{\theta 
							_k}} \right)}},\frac{{{L_H} + {\sigma _k} + {\kappa _\theta }{L_{\nabla 
							f}}}}{{1 - {\theta _k} 
					- {\zeta 
						_1} - {\zeta _2}}}} \right\}.\]
	\end{proof}
	%%===========================================================
	%%===========================================================
	%%===========================================================
	\textbf{Proof of Lemma {\ref{Newton:SARC:Lemma:bound_sigmal}}}
	\begin{proof}	
		We assume that ${\sigma _k} \ge \frac{9}{2}{L_H}$ and ${\sigma _k} \ge \frac{{\kappa 
				_H^2}}{{({\varepsilon _{\nabla f}} - {\varepsilon 
					_g})}}$, which are used for Lemma \ref{Newton:SARC:Lemma:UpperboundOfP0Ps} and 
		Lemma \ref{Newton:SARC:Lemma:Upperbound_p-h}. 	
		By setting $\frac{{\mathbb{I}\left( {\left| {{\mathcal{S}_h}} \right| < n} 
		\right)}}{{\left| 
				{{\mathcal{S}_h}} 
				\right|}}{H_1} \le {\varepsilon _g}$, and $\frac{{\mathbb{I}\left( {\left| 
					{{\mathcal{S}_h}} 
					\right| < n} 
				\right)}}{{\left| 
				{{\mathcal{S}_h}} \right|}}{H_2} \le {\varepsilon _B}$,  we consider 
		two cases:
		
		For the case $\nabla f\left( {{x_k}} \right) \ge {\varepsilon _{\nabla f}}$: Firstly, we 
		consider $g(x_k)$. Through adding and subtracting the term $\nabla f( x_k ) $, we 
		have the lower bound of $\|g(x)\|$,
		\begin{align*}
		\| {g( x_k )} \| &= \| {g( x_k ) -	 \nabla f( x_k ) + \nabla f( x_k )} 
		\|\nonumber\\
		&\ge \| {\nabla f( x_k )} \| -  \| {g( x_k 
			) - \nabla f( x_k )} \|\nonumber\\
		%\label{Newton:SARC:Lemma:Upperbound_mx-msx:inequality1}
		&\ge {\varepsilon _{\nabla f}} - {\varepsilon_g},
		\end{align*}
		where the last inequality is based on the approximation of $\nabla f(x_k)$ in Assumption 
		\ref{Newton:Assumption:approximation}. Secondly, because of 
		\[{\sigma _k} \ge \frac{{\kappa _H^2}}{{({\varepsilon _{\nabla f}} - {\varepsilon _g})}} 
		\ge 
		\frac{{{{\left\| {B({x_k})} \right\|}^2}}}{{\left\| {g({x_k})} \right\|}} \Rightarrow 
		\frac{{\left\| {B({x_k})} \right\|}}{{{\sigma _k}}} \le \sqrt {\frac{{\left\| {g({x_k})} 
					\right\|}}{{{\sigma _k}}}}, \]
		Combing with the upper bound of $s_k$ in Lemma \ref{Newton:SARC:Lemma:UpperboundOfP0Ps}, we 
		have 
		$\left\| 
		{{s_k}} \right\| \le \frac{{11}}{4}\sqrt {\left\| {g({x_k})} \right\|/{\sigma _k}}$.	
		Finally, based on equality (\ref{Newton:SARC:1-p}), Lemma 
		\ref{Newton:SARC:Lemma:UpperboundOfP0Ps} and Lemma \ref{Newton:SARC:Lemma:Upperbound_p-h}, 
		we 
		have
		\begin{align*}
		1 - \rho  
		\le& 10\frac{{4{\varepsilon _g}\left\| {{s_k}} \right\| + 3{\varepsilon _B}{{\left\| 
		{{s_k}} 
						\right\|}^2} + 2{\varepsilon _h}\|s_k\|^2}}{{\left\| {{g_k}} 
				\right\|\text{min} 
				\left\{ 
				{{\left\| {{g_k}} \right\|}/\left\| {{B_k}} \right\|,\sqrt {\left\| {{g_k}} 
						\right\|/{\sigma _k}} } \right\}}}
		\end{align*}
		In order to ensure that there exist a lower bound of $\sigma$ such that satisfying 
		$1-\rho\le 
		1-\eta$. Combing with the upper bound of $s_k$ in Lemma 
		\ref{Newton:SARC:Lemma:UpperboundOfP0Ps}, if ${\varepsilon _g} = 
		\frac{1}{{220}}\left( {1 - \eta } \right)({\varepsilon _{\nabla f}} - {\varepsilon 
			_g}),{\sigma _k} \ge \frac{1}{{{\varepsilon _{\nabla f}} - {\varepsilon 
					_g}}}\frac{{{{\left( {304\left( {3{\varepsilon _B} + 2{\varepsilon _h}} 
					\right)} 
						\right)}^2}}}{{\left( {1 - \eta } \right)}}$, we have
		\begin{align*}
		10\frac{{4{\varepsilon _g}\left\| {{s_k}} \right\| + \left( {3{\varepsilon _B} + 
					2{\varepsilon _h}} \right){{\left\| {{s_k}} \right\|}^2}}}{{\left\| {{g_k}} 
				\right\|{\rm{min}}\left\{ {\left\| {{g_k}} \right\|/\left\| {{B_k}} \right\|,\sqrt 
					{\left\| 
						{{g_k}} \right\|/{\sigma _k}} } \right\}}} \le& 10\frac{{11{\varepsilon 
					_g}\sqrt {\left\| 
					{g({x_k})} \right\|/{\sigma _k}}  + \left( {3{\varepsilon _B} + 2{\varepsilon 
					_h}} 
				\right){{\left( {\frac{{11}}{4}\sqrt {\left\| {g({x_k})} \right\|/{\sigma _k}} } 
						\right)}^2}}}{{\left\| {{g_k}} \right\|\sqrt {\left\| {{g_k}} 
						\right\|/{\sigma 
						_k}} }}\\
		\le& \frac{{110{\varepsilon _g}}}{{\left\| {{g_k}} \right\|}} + \frac{{76\left( 
				{3{\varepsilon _B} + 2{\varepsilon _h}} \right)\sqrt {\left\| {g({x_k})} 
					\right\|/{\sigma 
						_k}} }}{{\left\| {{g_k}} \right\|}} = \frac{{110{\varepsilon _g}}}{{\left\| 
				{{g_k}} 
				\right\|}} + \frac{{76\left( {3{\varepsilon _B} + 2{\varepsilon _h}} 
				\right)}}{{\sqrt 
				{\left\| {{g_k}} \right\|{\sigma _k}} }}\\
		\le& \frac{{110{\varepsilon _g}}}{{({\varepsilon _{\nabla f}} - {\varepsilon _g})}} + 
		\frac{{76\left( {3{\varepsilon _B} + 2{\varepsilon _h}} \right)}}{{\sqrt {({\varepsilon 
						_{\nabla f}} - {\varepsilon _g}){\sigma _k}} }} \le \frac{1}{2}\left( {1 - 
						\eta 
		} \right).
		\end{align*}
		Thus, we can see that if ${\sigma _{{\rm{max1}}}} = \frac{1}{{({\varepsilon _{\nabla f}} - 
				{\varepsilon _g})}}{r_2}\left\{ {\kappa _H^2,\frac{{{{\left( {304\left( 
				{3{\varepsilon 
										_B} + 
									2{\varepsilon _h}} \right)} \right)}^2}}}{{\left( {1 - \eta } 
					\right)}},\frac{9}{2}({\varepsilon _{\nabla f}} - {\varepsilon _g}){L_H}} 
		\right\}$, $r_2>1$, 
		we can obtain that 
		$1 - \rho  \le1-\eta$.

		For the case $ \lambda_\text{min} ({\nabla^2f(x )} ) \le - {\varepsilon_H }$: Firstly, 
		based on 
		the  Rayleigh quotient 
		\cite{conn2000trust} that if ${H\left( x \right)}$ is symmetric and the 
		vector $s \ne 0 $, 
		%	that is
		%	\begin{align}\label{SFSNewton:SARC:Lemma:UpperboundRayleigh_inequality1}
		%	{\lambda _{\text{min}  }}\left( {H\left( x \right)} \right) \le 
		%	\frac{{s_{}^TH\left( x \right)s}}{{{{\left\| s \right\|}^2}}} \le 
		%	{\lambda _{\text{max}  }}\left( {H\left( x \right)} \right),
		%	\end{align}
		then, we have
		\begin{align*}
		\frac{{s_k^TB\left( {{x_k}} \right){s_k}}}{{{{\left\| {{s_k}} 
						\right\|}^2}}} &= \frac{{s_k^T\left( {B\left( {{x_k}} 
					\right) - 
					{\nabla ^2}f\left( {{x_k}} \right) + {\nabla ^2}f\left( {{x_k}} \right)} 
				\right){s_k}}}{{{{\left\| 
						{{s_k}} \right\|}^2}}}\\
		&\mathop  \le \limits^{\scriptsize \textcircled{\tiny{1}}} 
		\frac{{s_k^T{\nabla ^2}f\left( {{x_k}} \right){s_k}}}{{{{\left\| {{s_k}} \right\|}^2}}} + 
		\left\| 
		{{\nabla ^2}f\left( {{x_k}} \right) - B\left( {{x_k}} \right)} \right\|\\
		&\mathop  \le \limits^{\scriptsize \textcircled{\tiny{2}}} {\lambda 
			_{\text{min}  }}\left( {{\nabla ^2}f\left( {{x_k}} \right)} \right) + {\varepsilon _B}
		\le  - 
		{\varepsilon _H } + {\varepsilon _B},
		\end{align*}
		where ${\scriptsize \textcircled{\tiny{1}}}$ is based on the triangle 
		inequality, 
		${\scriptsize \textcircled{\tiny{2}}}$ follows from the  Rayleigh quotient 
		\cite{conn2000trust}  and 
		approximation Assumption (\ref{Newton:Assumption:approximation-gB}). Secondly, based on 
		$s_k^T{B_k}{s_k} + {\sigma _k}{\left\| {{s_k}} \right\|^3} \ge 0$, we have ${\sigma 
		_k}\left\| 
		{{s_k}} \right\| \ge  - \frac{{s_k^T{B_k}{s_k}}}{{{{\left\| {{s_k}} \right\|}^2}}} \ge 
		\left( 
		{{\varepsilon _H} - {\varepsilon _B}} \right)$. Thirdly, based on equality 
		(\ref{Newton:SARC:1-p}), and  Lemma 
		\ref{Newton:SARC:Lemma:UpperboundOfP0Ps}, Lemma 
		\ref{Newton:SARC:Lemma:Upperbound_p-h}, we have
		\begin{align*}
		1 - \rho  
		%	\le&
		%	\frac{{\frac{3}{2}\left( {\frac{{\left( {\left| {{{\cal S}_h}} \right| < n} 
		%\right)}}{{\left| 
		%						{{{\cal S}_h}} \right|}}{H_2} + {\varepsilon _B}} \right){{\left\| 
		%{{s_k}} 
		%					\right\|}^2} + \left( 
		%			{\frac{3}{2}{L_H} - \frac{1}{3}{\sigma _k}} \right){{\left\| {{s_k}} 
		%\right\|}^3} + 
		%			2{\varepsilon 
		%				_h}\|s_k\|^2}}{{\frac{{{\sigma _k}}}{6}{{\left\| {{s_k}} \right\|}^3}}}\\
		\le& \frac{{2{\varepsilon _h}\|s_k\|^2}}{{\frac{{{\sigma _k}}}{6}{{\left\| {{s_k}} 
						\right\|}^3}}} + \frac{{3{\varepsilon _B}{{\left\| {{s_k}} 
						\right\|}^2}}}{{\frac{{{\sigma _k}}}{6}{{\left\| {{s_k}} \right\|}^3}}} 
		=\frac{{12{\varepsilon 
					_h}}}{{\sigma _k{{\left\| {{s_k}} \right\|}}}} + 18{\varepsilon 
			_B}\frac{1}{{{\sigma _k}\left\| 
				{{s_k}} \right\|}}\\
		\le& \frac{{12{\varepsilon _h}}}{{{{\left( {{\varepsilon _H} - {\varepsilon _B}} 
		\right)}}}} 
		+ 
		18{\varepsilon _B}\frac{1}{{{\varepsilon _H} - {\varepsilon _B}}},
		\end{align*}
		In order to have the lower bound radius $\Delta_k$ such that $1 - \rho  \le1-\eta$, we 
		consider 
		the parameters' setting:
		\begin{itemize}
			\item For the first term, then we define, 
			${\varepsilon _B} = \frac{1}{{36}}\left( {1 - \eta } \right)\left( {{\varepsilon _H} - 
				{\varepsilon _B}} \right)$.
			\item For the second term, we define ${\varepsilon _h} \le \frac{1}{{24}}{\left( 
				{{\varepsilon 
						_H} - {\varepsilon _B}} \right)}\left( {1 - \eta } \right)$,	
		\end{itemize}
		Thus, we can see that if ${\sigma _{\text{max2}}} = \frac{9}{2}{r_2}{L_H}$, $r_2>1$, we can 
		obtain 
		that 
		$1 - \rho  \le1-\eta$. 
		
		All in all, there is a large $\sigma$ such that lead to the successful iteration, that is 
		\[{\sigma _{\text{max}}} = \max \left\{ {{\sigma _\text{max1}},{\sigma _\text{max2}}} 
		\right\}.\]
	\end{proof}
	%%%==========================
	%%%===================================
	%%%====================================
	\textbf{Proof of Theorem \ref{Newton:SARC:theorem:Iteration}}
	%\begin{theorem}\label{Newton:SARC:theorem:Iteration}
	%	In Algorithm \ref{Newton:SARC:Algorithm}, suppose the Assumption 
	%	\ref{Newton:Assumption:function}-
	%	\ref{Newton:Assumption:Bound:Variance} hold, let $| {{\mathcal{S}_h}} | = \min 
	%	\{ n,\max \{ H_1/\varepsilon_g,H_2/\varepsilon_B 
	%	\} \}$, $\{ {f( {{x_k}} )} \}$ is 
	%	bounded below by $f_{\text{low}}$, the number of successful iterations ${T_\text{suc}}$  is 
	%	no large than
	%	\[{\kappa _5}\max \{ {{{({\varepsilon _{\nabla f}} - {\varepsilon _g})}^{ - 2}},{{( 
	%				{{\varepsilon _H} - {\varepsilon _B}} )}^{ - 3}}} \},\]
	%	where ${\kappa _5} = ( {f( {{x_0}} ) - {f_{{\rm{low}}}}} )\max \{ 
	%	{5/( \eta \kappa _4^{ - 1/2} ),6{\sigma _{{\rm{max2}}}}/\eta } \}$, 
	%	$\kappa_4$  is defined in (\ref{Newton:SARC:Lemma:bound_sigmal-2}). If the conditions 
	%	(\ref{Newton:SARC:Assumption_s1})-(\ref{Newton:SARC:Assumption_s3}) are satisfied, then the 
	%	number of successful 	iterations 
	%	$T_\text{unsuc}$ is at 
	%	most 
	%	\[{\kappa _6}\max \{ {{{({\varepsilon _{\nabla f}} - {\varepsilon _g})}^{ - 3/2}},{{( 
	%				{{\varepsilon _H} - {\varepsilon _B}} )}^{ - 3}}} \},\]
	%	where ${\kappa _6} = ( {f( {{x_0}} ) - {f_{{\text{low}}}}} )\max \{ 
	%	{6\kappa _s^{3/2}/( {\eta {\sigma _{{\text{min}}}}} ),6{\sigma _{{\text{max2}}}}/\eta } 
	%	\}$, $\kappa_s$ is defined (\ref{Newton:SARC:Lemma:UpperboundOfP0Ps-2-1}).
	%\end{theorem}
	\begin{proof}
		We consider two kinds of iteration complexity:
		\begin{itemize}
			\item 
			\begin{itemize}
				\item For the case of $\nabla f\left( {{x_k}} \right) \ge {\varepsilon _{\nabla 
				f}}$, 
				based on Lemma 
				\ref{Newton:SARC:Lemma:UpperboundOfP0Ps}, if $k$-th iteration is 
				successful, we obtain
				\begin{align*}
				f\left( {{x_k}} \right) - f\left( {{x_{k + 1}}} \right) \ge& {\eta}\left( {m\left( 
				0 
					\right) - {m_k}\left( s_k \right)} \right)\\
				\ge& \frac{1}{10}{\eta}({\varepsilon _{\nabla f}} - {\varepsilon _g})\sqrt 
				{\frac{{{\varepsilon _{\nabla f}} - {\varepsilon _g}}}{{{\sigma _{{\rm{max1}}}}}}}  
				= 
				\frac{1}{10}{\eta}{({\varepsilon _{\nabla f}} - {\varepsilon _g})^2}\kappa _4^{ - 
				1/2},
				\end{align*}
				where inequality  follow from 
				the inequality ${{\scriptsize \textcircled{\tiny{5}}}}$ is from $\left\| 
				{ f\left( {{x_k}} \right)} \right\| \ge {\varepsilon _{\nabla 
						f}}-{\varepsilon _{g}}$ and $\sigma_{\text{max1}}$ in 
				(\ref{Newton:SARC:Lemma:bound_sigmal-1}). Let $T_3$ is the number of  successful 
				iteration, we obtain
				\begin{align*}
				f\left( {{x_0}} \right) - f_{\text{low}} 
				%	&\ge \sum\limits_{k = 
				%		0,k \in S}^{} {f\left( {{x_k}} \right) - f\left( {{x_{k + 1}}} 
				%		\right)} 
				%	\\
				%	&
				\ge 
				\frac{1}{5}T_3{\eta}{({\varepsilon _{\nabla f}} - 
					{\varepsilon 
						_g})^2}\kappa_4^{ - 1/2}.
				\end{align*}
				\item For the case of ${\lambda _{\text{min} }}( \nabla^2f(x_k) ) 
				\le  - 	{\varepsilon _H}$, based on Lemma 
				\ref{Newton:SARC:Lemma:UpperboundOfP0Ps-2}, 
				if 
				$k$-th 
				iteration is 
				successful, we have
				\begin{align*}
				f\left( {{x_k}} \right) - f\left( {{x_k} + {s_k}} \right) &
				\ge  \frac{{{\sigma _k}}}{6}\eta{\left\| {{s_k}} 
					\right\|^3} = \frac{1}{{6\sigma _k^2}}\eta\sigma _k^3{\left\| {{s_k}} 
					\right\|^3}\\
				& \ge - \frac{1}{{6{\sigma 
							_\text{max2}}}}\eta{\left( {\frac{{s_k^T{B_k}{s_k}}}{{{{\left\| {{s_k}} 
										\right\|}^2}}}} 
					\right)^3}\\
				& \ge  \frac{1}{{6{\sigma 
							_\text{max2}}}}\eta{\left( {{\varepsilon _H} - {\varepsilon _B}} 
							\right)^3},
				\end{align*}
				where  inequalities are based on  the condition of $s_k$ that 
				satisfies 
				(\ref{Newton:SARC:Assumption_s2}),  that is 
				$
				{\sigma _k}\left\| {{s_k}} \right\| \ge  - \frac{{s_k^TB\left( {{x_k}} 
						\right){s_k}}}{{{{\left\| {{s_k}} \right\|}^2}}},
				$ and $\sigma_{\text{max2}}$ in (\ref{Newton:SARC:Lemma:bound_sigmal-1}). Let $T_4$ 
				be 
				the number of successful iteration, then we obtain	
				\[f( {{x_0}} ) - f_\text{low} \ge \frac{1}{{6{\sigma _\text{max2}}}}T_4\eta {\left( 
					{{\varepsilon _H} - 
						{\varepsilon _B}} \right)^3}.\]	
			\end{itemize}
			Thus the the total number of success iteration is \[{T_3} + {T_4} = {\kappa _5}\max 
			\left\{ 
			{{{({\varepsilon _{\nabla f}} - {\varepsilon _g})}^{ - 2}},{{\left( {{\varepsilon _H} - 
							{\varepsilon _B}} \right)}^{ - 3}}} \right\}\]
			where ${\kappa _5} = \left( {f\left( {{x_0}} \right) - {f_{{\rm{low}}}}} \right)\max 
			\left\{ 
			{5/\left( {\eta \kappa _4^{ - 1/2}} \right),6{\sigma _{{\rm{max2}}}}/\eta } \right\}$
			\item  		For the case that condition 
			(\ref{Newton:SARC:Assumption_s1})-(\ref{Newton:SARC:Assumption_s3}) and $
			g(x_k)>(\varepsilon_{\nabla f}-\varepsilon_g)$ hold, we have
			\begin{align*}
			f\left( {{x_k}} \right) - {f}\left( {{x_k} + {s_k}} \right) \ge 
			\frac{{{\sigma 
						_k}}}{6}{\left\| {{s_k}} \right\|^3}
			\ge\frac{{{\sigma _{\text{min} }}}\kappa_s^{-3/2}}{6}{\left\| {g\left( {{x_k}} \right)} 
				\right\|^{3/2}}
			\ge \frac{{{\sigma _{\text{min} }}}\kappa_s^{-3/2}}{6}{\left( {{\varepsilon _{\nabla 
			f}} - 
					{\varepsilon _g}} 
				\right)^{3/2}}
			\end{align*}
			where  inequality are based on Lemma \ref{Newton:SARC:Lemma:UpperboundOfP0Ps-2} and 
			$\sigma_{{\text{min}}}$ in (\ref{Newton:SARC:Lemma:bound_sigmal-1}). Let $T_5$ is the 
			number 
			of the successful 
			iterations, we obtain	
			\begin{align*}
			f({x_0}) - {f_\text{low}} \ge \frac{{{\sigma _{\text{min} }}}\kappa_s^{-3/2}}{6}T_5\eta 
			{\left( 
				{{\varepsilon 
						_{\nabla 
							f}} - {\varepsilon _g}} 
				\right)^{3/2}}.
			\end{align*}
			The total number of success iteration  for such case is
			\[{T_3} + {T_5} = {\kappa _6}\max \left\{ {{{({\varepsilon _{\nabla f}} - {\varepsilon 
							_g})}^{ 
						- 3/2}},{{\left( {{\varepsilon _H} - {\varepsilon _B}} \right)}^{ - 3}}} 
			\right\},\]
			where ${\kappa _6} = \left( {f\left( {{x_0}} \right) - {f_{{\rm{low}}}}} \right)\max 
			\left\{ {6\kappa _s^{3/2}/\left( {\eta {\sigma _{{\rm{min}}}}} \right),6{\sigma 
					_{{\rm{max2}}}}/\eta } \right\}$.
		\end{itemize}
	\end{proof}

	%%%%%%%%%%%%%%%%%%%%%%%%%%%%%%%%%%%%%%%%%%%%%%%%%%%%%%%%%%%%%
	%%%%%%%%%%%%%%%%%%%%%%%%%%%%%%%%%%%%%%%%%%%%%%%%%%%%%%%%%%%
	%%%%%%%%%%%%%%%%%%%%%%%%%%%%%%%%%%%%%%%%%%%%%%%%%%%%%%%%%%
	%\input{Section-Appendix-STR}
	%\input{Section-Appendix-SARC}
	
\end{document}